\documentstyle[12pt]{article}
\title{Actions of $GL_q(2,C)$ on  $C(1,3)$ and its
four dimensional representations}
\author{V.K. Kharchenko, Jaime Keller and 
S. Rodr\'{\i}guez-Romo $^{\ast}$\\
Centro de Investigaciones Te\'oricas,\\ 
Universidad Nacional Aut\'onoma de M\'exico, Campus Cuautitl\'an \\
Apdo. Postal 95, Cuautitl\'an Izcalli, 
Edo. de M\'exico, 54768 M\'exico}         
\date{}  
\begin{document}   
\maketitle 
\renewcommand{\thefootnote}{\fnsymbol{footnote}}
\setcounter{footnote}{-1}
\footnote{$\hspace*{-6mm}^{\ast}$
e-mail: suemi$@$servidor.unam.mx }
\renewcommand{\thefootnote}{\arabic{footnote}}\
\renewcommand{\ref}{}
\baselineskip0.6cm
{\small {\bf Abstract.} A complete classification is given 
of all inner actions
on the Clifford algebra $C(1,3)$
defined by representations of the quantum group $GL_q(2,C),$ \ $q^m\neq 1$
 with nonzero perturbations. As a consequence of this classification
 it is shown that the space of invariants of every $GL_q(2,C)$-action 
of this type, which is not an action of $SL_q(2,C),$ is generated by 1 and 
the value of the quantum determinant for the given representation.} 
\newtheorem{theorem}{Theorem}
\newtheorem{corollary}{Corollary}
\newtheorem{lemma}{Lemma}  
\section{Introduction}
The development of the quantum group theory and its applications
in mathematical physics (more exactly an interpretation of
some results using physical intuitions) give us the
hope  that the notions of quantum actions, quantum invariants 
and quantum symmetries could play an important role 
in the ``quantum" mathematics as the classical notions
of invariants and symmetry do in classical mathematics 
and theoretical physics.
In this paper we study in details actions and invariants
of the quantum group $GL_q(2,C)$ on the space-time Clifford 
algebra $C(1,3),$ which is 
generated by vectors $\gamma _{\mu }, \mu =0,1,2,3$ 
with relations defined by the form
$g_{\mu \nu}=(1,-1,-1,-1).$ This algebra
has 16 matrix units 
\begin{equation}
\vbox{\halign{# \hfil &\ # \hfil \cr
$e_{11}=(1+\gamma_0)(1+i\gamma_{12})/4,$&
$e_{21}=(1+\gamma_0)(i\gamma_2-\gamma_1)\gamma_3/4,$\cr
$e_{31}=(1-\gamma_0)(1+i\gamma_{12})\gamma_3/4,$ & 
$e_{41}=(1-\gamma_0)(\gamma_1-i\gamma_2)/4,\nonumber $\cr
$e_{12}=(1+\gamma_0)(\gamma_1+i\gamma_2)\gamma_3/4,$ & 
$e_{22}=(1+\gamma_0)(1-i\gamma_{12})/4,\nonumber $\cr
$e_{32}=(1-\gamma_0)(\gamma_1+i\gamma_2)/4,$ &
$e_{42}=(1-\gamma_0)(i\gamma_{12}-1)\gamma_3/4,\nonumber $\cr
$e_{13}=-(1+\gamma_0)(1+i\gamma_{12})\gamma_3/4,$ & 
$e_{23}=-(1+\gamma_0)(\gamma_1-i\gamma_2)/4,\nonumber $\cr
$e_{33}=(1-\gamma_0)(1+i\gamma_{12})/4,$ &
$e_{43}=(1-\gamma_0)(-\gamma_1+i\gamma_2)\gamma_3/4,\nonumber $\cr
$e_{14}=-(1+\gamma_0)(\gamma_1+i\gamma_2)/4,$ &
$e_{24}=(1+\gamma_0)(1-i\gamma_{12})\gamma_3/4,\nonumber $\cr
$e_{34}=(1-\gamma_0)(\gamma_1+i\gamma_2)\gamma_3/4,$ & 
$e_{44}=(1-\gamma_0)(1-i\gamma_{12})/4.\nonumber $\cr }}
\label{matrix}
\end{equation}
and therefore it is
abstractly isomorphic to the algebra of 4 by 4 complex matrices.

The quantum group ${\it GL}_q(2,C)$ is generated 
by four, so-called, $q$-spinors; therefore the 
classification of representations of the $q$-spinor in $C(1,3)$
given in Theorem \ref{t2} 
produces a basic tool for the investigation of the inner actions defined by 
representations of the algebraic structure of ${\it GL}_q(2,C)$. 
In the fourth section we
find a necessary and sufficient condition for different 
representations to 
define equal or equivalent actions. Using this result in Sections 6, 7 we 
have obtained a complete classification  
of all inner actions
defined by representations 
with nonzero ``perturbations" 
(i.e. when two main $q$-spinors do not commute).
Results are summarized in the Table 
given in Section 8. 
As a consequence of this classification
we have shown that the space of quantum invariants of 
every $GL_q(2,C)$-action 
of this type, which is not an action of $SL_q(2,C),$ is generated by 1 and 
the value of the quantum determinant. Roughly speaking, it means 
that in this case the quantum determinants are the only quantum invariants.

We believe that the results of this classification,  
besides possible interpretations
with physical intuition, could be usefull as a starting material 
for investigations of actions of more complicated quantum groups on another 
natural classical objects.

\section{Preliminary notions} 
The quantum group $GL_q(2,C)$ is a Hopf algebra whose algebraic structure
is generated by five elements $a_{11}, a_{12}, a_{21}, a_{22}, d^{-1}$
which can be presented by the following diagram: 
\begin{equation} 
\begin{picture}(80,80)
\put(15.29,2.39){\makebox(0,0)[cc]{$a_{21}$}}
\put(60.0,2.59){\makebox(0,0)[cc]{$a_{22}$}}
\put(60.0,45.46){\makebox(0,0)[cc]{$a_{12}$}}
\put(15.29,45.46){\makebox(0,0)[cc]{$a_{11}$}}
\put(22,45.46){\vector(1,0){28}} 
\put(22,2.39){\vector(1,0){28}}
\put(15.29,38){\vector(0,-1){28}}
\put(58.32,38){\vector(0,-1){28}}
\put(20.0,8.0){\line(1,1){30}}
\put(25.0,33.0){$\cdot$}
\put(30.0,28.0){$\cdot$}
\put(35.0,23.0){$\cdot$}
\put(40.0,18.0){$\cdot$}
\put(45.0,13.0){$\cdot$}
\put(50.0,8.0){$\cdot$}
\end{picture}
\label{p1}
\end{equation}
where by the arrows $x\rightarrow y$ 
are denoted the so-called $q$-spinors $xy=qyx$, 
by the straight line commuting elements
and by dots, elements with nontrivial commutator 
$[a_{11}a_{22}]=(q-q^{-1})a_{12}a_{21}$
(for another quantum deformations of $GL_n$ see 
\cite{Art},\cite{Dem} \cite{Dip}, \cite{Wor}).
The comultiplication and the counit are defined as follows:
$
\Delta(a_{ij})=\sum_{k=1}^{k=2}a_{ik}\otimes a_{kj},$\ 
$\ \varepsilon (a_{ij})=\delta ^j_i.$
The antipode is given by the following formula
\begin{equation}
S\left(\matrix{a_{11}&a_{12}\cr a_{21}&a_{22}\cr}\right)=
d^{-1}\left(\matrix{a_{22}&-q^{-1}a_{12}\cr -qa_{21}&a_{11}\cr}\right).
\label{ant}
\end{equation}

The quantum group $SL_q(2,C)$ is 
defined as  the factor-Hopf algebra of $GL_q(2,C)$
by the additional relation $d=1.$

An action of Hopf algebra $H$ on an algebra $A$ is characterized 
by the following two main formulae
\begin{equation}
(hg)\cdot v=h\cdot (g\cdot v),
\label{mod}
\end{equation}
\begin{equation}
h\cdot vw=\sum (h_{(1)}\cdot v)(h_{(2)}\cdot w),
\label{act}
\end{equation}
where $h, g\in H, \ v, w\in A$ and $\Delta (h)=\sum h_{(1)}\otimes h_{(2)}$
(see details in \cite{Coh}, \cite{Shn}).

The first formula shows that 
the action will be defined if it is defined for generators $a_{ij}$ of
$GL_q(2,C)$, while the second is showing that for a 
definition of an action it 
is enough to set the action of $a_{ij}$ on the generators $\gamma_0$, 
$\gamma_1$, $\gamma_2$, $\gamma_3$ of the Clifford algebra  ${\it C}(1,3)$. 
Thus, an action is set by formulas of the following type
\begin{equation}
a_{ij}\cdot \gamma_k=f_{ijk}(\gamma_0, \gamma_1, \gamma_2, \gamma_3)
\label{fact}
\end{equation}
where $f_{ijk}$ are some noncommutative polynomials in four variables. 

Two actions $*$ and $\cdot $ of a Hopf algebra 
$H$ on an algebra $A$
are called {\it equivalent} if 
$
h*v=
\left( h\cdot v^{\zeta^{-1}}\right)^{\zeta},
$
where $\zeta$ is an automorphism of the algebra $A.$
By Skolem---Noether theorem every automorphism of $C(1,3)$
is given by a conjugation $w^{\zeta}=uwu^{-1}.$ Therefore equivalence of
actions can be presented by the following formula
\begin{equation}
a_{ij}*(uwu^{-1})=u(a_{ij}\cdot w)u^{-1}.
\label{fequiv}
\end{equation}
This formula shows that for an action $*$ to be equivalent   
to the action $\cdot $ there exists a presentation
\begin{equation}
a_{ij}*\gamma'_k=f_{ijk}(\gamma'_0, \gamma'_1, \gamma'_2, \gamma'_3)
\end{equation}
with the same polynomials $f_{ijk}$ and a
system of generators
$\gamma'_0$, $\gamma'_1$, $\gamma'_2$, $\gamma'_3$ 
with the same relations.

Using Skolem---Noether theorem 
(see also \cite{Ko}, or \cite{Mo} chapter 6.2.) it is easy to show that for 
every action there exists an invertible 2 by 2 matrix $M$ over $C(1,3)$
such that 
\begin{equation}
\left(\matrix{a_{11}\cdot v &a_{12}\cdot v \cr
              a_{21}\cdot v &a_{22}\cdot v \cr}\right) =
M\left(\matrix{v&0\cr 0&v\cr }\right)M^{-1}.
\label{sko}
\end{equation}

Inversely, if $M=\left(\matrix{m_{11}&m_{12}\cr m_{21}&m_{22}\cr }\right)$
is a 2 by 2 matrix with the additional condition 
that $a_{ij}\rightarrow A_{ij} =
m_{ij}$ represents a homomorphism of the algebraic 
structure of $GL_q(2,C)$ to
$C(1,3),$ then $M$ is invertible and (\ref{sko}) defines an action of
$GL_q(2,C)$ on $C(1,3)$ (called an {\it inner} action). 

For a given representation $a_{ij}\rightarrow A_{ij}$ we denote by
$\Re $ an {\it operator algebra} i.e. a 
subalgebra of $C(1,3)$ generated by $A_{ij}.$ 
Recall that the algebra of invariants of an action is defined in the
following way
\begin{equation}
Inv=\{ v\in A | \forall h\in H \ \ \ h\cdot v=\varepsilon (h)v\} .
\label{inv}
\end{equation}

\begin{lemma}. If an action of $GL_q(2,C)$ is given by the formula
$(\ref{sko})$ then the algebra $Inv$ equals the centralizer of 
all components of $M.$ In particular the algebra $Inv$ for an inner 
action on $C(1,3)$ defined by a representation $a_{ij}\rightarrow A_{ij}$
equals the centralizer of $\Re $ in $C(1,3).$
\end{lemma}
{\it Proof.}  An element $v$ is an invariant if an only if
$a_{11}\cdot v=\varepsilon (a_{11})v=v,$\ 
$a_{12}\cdot v=\varepsilon (a_{12})v=0,$\
$a_{21}\cdot v=\varepsilon (a_{21})v=0,$\
$a_{22}\cdot v=\varepsilon (a_{22})v=v.$
In matrix form this is equivalent to
$diag(v,v)=M\, diag(v,v)M^{-1}.$ Therefore $diag(v,v)M=M\, diag(v,v),$ 
i.e. $vm_{ij}=m_{ij}v$ for all components $m_{ij}$ of $M.$
\hfill$\Box$
\section{$q$-spinor representations}
Let $(x,y)\in A^{2/0}_q$ be a $q$-spinor $xy=qyx$. If 
$x\rightarrow A$, $y\rightarrow B$ is its representation by $4\times 4$ 
matrices over complex numbers, then for every invertible $4\times 4$ matrix 
$u$ and nonzero number $\alpha,$ the map $x\rightarrow uAu^{-1}\alpha$, 
$y\rightarrow uBu^{-1}\alpha$ also defines a representation of the $q$-spinor. 
We consider this two representations as {\it equivalent} ones. 
Thus, under investigation of representations of a $q$-spinor, we can 
suppose that the matrix $A$ has a Jordan Normal form and one of it's 
eigenvalues is equal to 1 (if $A\neq 0)$.\

For a given matrix $A,$ a set $B(A)$ of all matrices $B,$ such that 
$x\rightarrow A$, $y\rightarrow B$ is a representation of the 
$q$-spinor, forms a linear space. Therefore, it is natural to represent the
space $B(A)$ by one of it's basis: $\{B_1, B_2, ...\}.$
\begin{theorem}.
Every representation of the $q$-spinor $(q^3, q^4\neq 1)$ by $4\times 4$
complex matrices, $x\rightarrow A$, $y\rightarrow B$, such that $A$ is
invertible and $B(A)^2\neq 0,$ is equivalent to one of the following 
representations
\label{t2}
\end{theorem}
\begin{equation}
1.\  A=diag(q^2, q, 1, 1), \;\;B_1=e_{12},\;B_2=e_{23},\;B_3=e_{24}\hfill 
\label{f31}
\end{equation}
\begin{equation}
2.\ A=diag(q^2, q, q, 1),\;\; B_1=e_{12},\;B_2=e_{13},\;B_3=e_{24},
\;B_4=e_{34}\hfill 
\label{f32}
\end{equation}
\begin{equation}
3. \  A=diag(q^2, q^2, q, 1),\;\; B_1=e_{13},
\;B_2=e_{23},\;B_3=e_{34}\hfill 
\label{f33}
\end{equation}
\begin{equation}
4. \ A=diag(q^3, q^2, q, 1),\;\;B_1=e_{12},\;B_2=e_{23},\;B_3=e_{34}\hfill 
\label{f34}
\end{equation}
\begin{equation}
5. \ A=diag(\alpha, q^2, q, 1),\;\;B_1=e_{23},\;B_2=e_{34},\ 
\alpha \neq 0, q^{-1}, 1, q, q^2, q^3\hfill  
\label{f35}
\end{equation}
\begin{equation}
6. \ A=diag(q^2, q^2, q, 1)+e_{12},\;\;B_1=e_{13},\;B_2=e_{34}\hfill 
\label{f36}
\end{equation}
\begin{equation}
7. \ A=diag(q^2, q, 1, 1)+e_{34},\;\;B_1=e_{24},\;B_2=e_{12}\hfill 
\label{f37}
\end{equation}
\begin{lemma}.
Let a matrix $A$ have the form $A$=$\alpha E+U$, where 
$\alpha\neq  0$, $U$ is a nilpotent $n\times n$-matrix, $E$ is the identity 
matrix, $E$=diag$(1,1,...,1)$. If $AB$=$qBA$, $q\neq 1$, then $B=0$.
\label{2}
\end{lemma}
{\it Proof}. If we denote $[x,y]_q$=$xy-qyx$ then $[A,B]_q=0$; 
i.e. 
\begin{equation}
[\alpha E+U,B]_q=\alpha(1-q)B+[U,B]_q.
\label{f40}
\end{equation}
This formula implies
$
[U,B]_q=\alpha(q-1)B.
$
Iterating this equality $m$ times we have
\begin{equation}
[U,[U,[...[U,B]_q...]_q]_q=\alpha^m(q-1)^mB. \label{f42}
\end{equation}
If $U^k=0$ and $m=2k$ then the left hand part of (\ref{f42}) is equal to zero. Thus
$\alpha^m(q-1)^mB$=$0$ and $B=0$. \hfill$\Box$\
\begin{lemma}
. If $x\rightarrow A$, $y\rightarrow B$ is a representation of 
the $q$-spinor and $C_1$, $C_2$ are matrices commuting with $A$, then 
$x\rightarrow A$, $y\rightarrow C_1BC_2$ is also a representation of the 
$q$-spinor.
\label{3}
\end{lemma}
{\it Proof}. $A\cdot C_1BC_2$=$C_1A\cdot BC_2$=$C_1qBA\cdot C_2$=
$qC_1BC_2\cdot A$.\hfill$\Box$\

\begin{lemma}. If $A$ is a diagonal matrix then $B(A)$ has a basis consisting of matrix
units.
\label{4}
\end{lemma}
{\it Proof}. If $B=||\beta_{ij}||\in B(A)$, then $\beta_{ij}e_{ij}$=
$e_{ii}Be_{jj}\in B(A)$ (since $e_{ii}A$=$Ae_{ii}$ for 
each diagonal matrix $A$ 
and we can use Lemma \ref{3}). Thus, for every nonzero 
$\beta_{ij}$ the matrix unit $e_{ij}$ belongs to $B(A)$.\hfill$\Box$

{\it Proof of Theorem} \ref{t2}. We can assume that the matrix $A$ has a 
Jordan Normal form and one of it's eigenvalues is equal to 1.\

Assume firstly that $A$ is a diagonal matrix; i.e. it has four simplest
Jordan Normal blocks: $A$=diag$(\alpha_1, \alpha_2, \alpha_3, \alpha_4)$,
$\alpha_4=1$. In this case, by Lemma \ref{4}, it is enough to describe all matrix 
units from $B(A)$. We have
\begin{equation}
Ae_{ij}=\alpha_ie_{ij}=\alpha_i\alpha^{-1}_j(\alpha_je_{ij})=
\alpha_i\alpha^{-1}_j(e_{ij}A).\label{f43}
\end{equation}
This formula shows that $e_{ij}\in B(A)$ if and only if 
$\alpha_i\alpha^{-1}_j=q$; i.e.
\begin{equation}
\alpha_i=q\alpha_j.\label{f44}
\end{equation}
If $B(A)^2\neq 0$, then there exist indices $i, j, k$ such that $e_{ij}$,
$e_{jk}\in B(A)$. A conjugation by matrix $T=E-e_{ii}-e_{jj}+e_{ij}+e_{ji}$
changes i'th and j'th entries in any diagonal matrix. Therefore we can suppose
that $i=2,j=3,k=4$; i.e. (\ref{f44}) 
implies $A$=diag$(\alpha, q^2, q, 1)$. If 
$\alpha\neq q^{-1}, 1,q, q^2, q^3$ then $B(A)$=$Ce_{23}+Ce_{34}$ and we have 
a representation (35). If $\alpha= q^{-1}, 1,q, q^2, q^3$ 
then we will obtain four possibilities which are equivalent to (\ref{f34}), 
(\ref{f31}),
(\ref{f32}), (\ref{f33}) and (\ref{f34}), respectively.\

Now, let us consider the cases such that $A$ has less than four blocks.\

By Lemma \ref{2}, $A$ cannot be a simplest 
Jordan Normal matrix; i.e. it has more 
than one block. This means that $A$ has a form
$
A=diag(a,b),
$
where $a,b$ are either invertible $2\times 2$ matrices in Jordan Normal form
or $a$ is an invertible simplest Normal Jordan $3\times 3$ matrix and $b$ is
a nonzero complex number (and therefore we can suppose that $b=1$). Let us
write
\begin{equation}
B=
\left(
\begin{array}{cc}
\alpha & \beta \\
\gamma & \delta
\end{array}
\right)\neq 0,
\label{f46}
\end{equation}
where $\alpha$, $\beta$, $\gamma$, $\delta$ are, respectively, either 
$2\times 2$ matrices or $\alpha$ is a $3\times 3$ matrix, $\delta$ is a 
complex number and $\beta$, $\gamma$, a column and a row, respectively. In 
both cases we have
\begin{equation}
[A,B]_q= 
\left(
\begin {array}{cc}
a\alpha-q\alpha a  & a\beta-q\beta b \\
b\gamma-q\gamma a  & b\delta-q\delta b
\end{array}
\right)=0.
\label{f47}
\end{equation}
This implies, in particular, that $x\rightarrow a$, $y\rightarrow \alpha$ and
$x\rightarrow b$, $y\rightarrow \delta$ are representations of 
$q$-spinors.\

Let us consider firstly the case when $a$ is a $3\times 3$ matrix. In this 
case, by Lemma \ref{2}, we have $\alpha=0$, $\delta=0$ and also $a\beta=q\beta$,
$\gamma=q\gamma a$ (recall that $b=1$). Therefore for a $3\times 3$ matrix
$\beta\gamma$ (this is the matrix product of a column by a row), we have
$
a(\beta\gamma)=q\beta\gamma=q^2(\beta\gamma)a.
$
Again by Lemma \ref{2} we obtain $\beta\gamma=0,$ therefore either $\beta=0$ or
$\gamma=0$.\

Let $\beta=0$. Then $\gamma\neq 0$ and the equality $\gamma=q\gamma a$ can
be written in the form $\gamma a$=$q^{-1}\gamma$. It means that $q^{-1}$ is 
an eigenvalue of $a$ and $\gamma$ is an eigenvector of $a$; i.e.

In this case $B(A)$=$Ce_{43}$ and $B(A)^2$=$0$.
$a=q^{-1}E+e_{12}+e_{23},$\ $ \gamma =\epsilon e_{43}.$
Let $\gamma=0,$ then $\beta\neq 0$ and the equality $a\beta=q\beta$ shows that
$q$ is the eigenvalue of $a$ and $\beta$ is an eigenvector i.e.
$a=qE+e_{12}+e_{23},$\ $\beta =\epsilon e_{14}$
and $B(A)$=$Ce_{14},$ so $B(A)^2=0$.

Consider now the case when $a, b, \alpha, \beta, \gamma, \delta$ are 
$2\times 2$ matrices.\

Let us start with the case when both matrices $a$ and $b$ have a simplest
Jordan Normal Form i.e.
$a=\epsilon E+e_{12},$\ $b=E+e_{34}$
By formula (\ref{f47}) and Lemma \ref{2} we have $\alpha=\delta=0$. If
\begin{equation}
B'=
\left(
\begin{array}{cc}
\alpha' & \beta' \\
\gamma' & \delta'
\end{array}
\right)
\label{f54}
\end{equation}
is another matrix from $B(A)$; then we also have $\alpha'=\delta'=0$ and by
(\ref{f47}) we have the relations
$
a\beta'=q\beta' b,\;\;\; b\gamma'=q\gamma'a.
$
This relations and (\ref{f47}) imply
$
a(\beta\gamma')=q\beta b\gamma'=q^2(\beta\gamma')a
$
and by Lemma \ref{2}, $\beta\gamma'=0$.\

In the same way
$
b(\gamma \beta')=q\gamma a\beta'=q^2(\gamma\beta')b
$
and by Lemma 2, $\gamma\beta'=0$. Now
\begin{equation}
BB'=
\left(
\begin{array}{cc}
0       & \beta \\
\gamma  & 0 
\end{array}
\right)
\left(
\begin{array}{cc}
0       & \beta'\\
\gamma'  & 0 
\end{array}
\right)=
\left(
\begin{array}{cc}
\beta\gamma' & 0 \\
0  & \gamma\beta' 
\end{array}
\right)=0,
\label{f58}
\end{equation}
therefore $B(A)^2$=$0$.

Suppose now that one of the matrix $a,b$ is a simplest Jordan matrix while 
another is a diagonal matrix. A conjugation by 
$T=
\left(
\begin{array}{cc}
0  & E\\
E  & 0
\end{array}
\right),$
where $E$ is the identity $2\times 2$ matrix, $E$=diag$(1,1)$, changes
$A$=diag$(a,b)$ to diag$(b,a)$, so we can suppose that
$
a=\epsilon E+e_{12},\;\;\;b=diag(\mu,1),
$
where  $\epsilon, \mu\neq 0$.

By (\ref{f47}) $x\rightarrow a$, $y\rightarrow \alpha$ 
is a representation of the
$q$-spinor and Lemma \ref{2} implies $\alpha=0$. 
In the same way $x\rightarrow b$,
$y\rightarrow \delta$ is a representation of the $q$-spinor. Therefore
if $\mu\neq q$, $q^{-1}$ then $\delta=0$ and we can repeat
word by word the proof beginning from formula (\ref{f54}) up to 
(\ref{f58}) in order to show that $B(A)^2=0$.\

If $\mu=q^{-1}$ we can multiply matrices $A$ and $B$ by $q$ and conjugate
them by the matrix
$
diag(1,1,
\left(
\begin{array}{cc}
0 & 1 \\
1 & 0
\end{array}
\right)).
$
We will obtain an equivalent representation with $\mu=q$. Thus, it is enough
to consider the case
\begin{equation}
a=\epsilon E+e_{12},\;\;\; b=diag(q,1)
\label{f61}
\end{equation}
and by Lemma \ref{2}
\begin{equation}
\alpha=0,\;\;\; \delta=ce_{12}, \;\;\;c\in \hbox{\bf C}.
\label{f62}
\end{equation}

For 
$\beta=
\left(
\begin{array}{cc}
\beta_{11} & \beta_{12}\\
\beta_{21} & \beta_{22}
\end{array}
\right),$ by (\ref{f47}), we have $\alpha \beta =q\beta b\;$,
which implies
\begin{equation}
(\epsilon-q^2)\beta_{11}=-\beta_{21}, \ \  
(\epsilon-q)\beta_{12}=-\beta_{22}
\label{f64}
\end{equation}
\begin{equation}
(\epsilon -q^2)\beta_{21}=0,         \ \  (\epsilon-q)\beta_{22}=0.
\label{f65}
\end{equation}
If $\epsilon$=$q^2$ then the first equality of (\ref{f64}) 
gives $\beta_{21}=0$, 
and if $\epsilon\neq q^2$ then the first equality of (\ref{f65}) gives 
$\beta_{21}=0$. Therefore $\beta_{21}=0$ in any case. In the same way 
$\beta_{22}=0$ and (\ref{f64}), (\ref{f65}) are equivalent to
\begin{equation}
(\epsilon-q^2)\beta_{11}=0 \ \ \  (\epsilon-q)\beta_{12}=0
\label{f66}
\end{equation}
\begin{equation}
\beta_{21}=0,          \ \ \  \beta_{22}=0.
\label{f67}
\end{equation}
Analogously for the matrix 
$\gamma=
\left(
\begin{array}{cc}
\gamma_{11} & \gamma_{12} \\
\gamma_{21} & \gamma_{22}
\end{array}
\right)$
we have $b\gamma$=$q\gamma a$; 
which is equivalent to
\begin{equation}
\gamma_{11}=0, \ \ \  (1-\epsilon)\gamma_{12}=0,
\label{f71}
\end{equation}
\begin{equation}
\gamma_{21}=0, \ \ \  (1-q\epsilon)\gamma_{22}=0.
\label{f72}
\end{equation}
Now if $\epsilon\neq q^{-1},1,q,q^2$ then by (\ref{f66}),(\ref{f67})
 and (\ref{f71}), (\ref{f72}) 
$\beta=\gamma=0$ and
the representation has the form
\begin{equation}
A=diag\left(
\left(
\begin{array}{cc}
\epsilon & 1 \\
0	 & \epsilon\\
\end{array}
\right), q, 1\right),\;\;\; B_1=e_{34}.
\label{f73}
\end{equation}
This means $B(A)=Ce_{34}$ and $B(A)^2=0$.\

Finally, let us consider four last possibilities.\

1. $\epsilon=q^{-1}$. By (\ref{f66}) and  (\ref{f67}) we have 
$\beta=0$ and by (\ref{f71}) and (\ref{f72}), 
$\gamma=ce_{22}$. Together with (\ref{f61}) and (\ref{f62}) it follows
\begin{equation}
A=diag\left(
\left(
\begin{array}{cc}
q^{-1} & 1 \\
0      & q^{-1}\\
\end{array}
\right),q,1\right), \;\;\;B_1=e_{42} , \;\;\;B_2=e_{34}.
\label{f74}
\end{equation}
If we multiply $A$ by $q$ and conjugate it by $T$=diag$(1, q^{-1},1,1)$ we 
will obtain an equivalent representation $A$=$diag(1,1,q^2,q)+e_{12}$, 
$B_1$=$e_{42}$, $B_{2}$=$e_{34}$. Using conjugations by matrices 
$E-e_{ii}-e_{jj}+e_{ij}+e_{ji}$ we can change indices with the help of
permutation $1\rightarrow 3$, $2\rightarrow 4$, $3\rightarrow 1$,
$4\rightarrow 2$. In this case $e_{42}\rightarrow e_{24}$, $e_{34}\rightarrow
e_{12}$ and we have the representation (\ref{f37}).\

2. $\epsilon=1$. By (\ref{f66}) and (\ref{f67})
we again have $\beta=0$ and by (\ref{f71}) and (\ref{f72}), 
$\gamma=Ce_{12}$. Now relations (\ref{f61}) and (\ref{f62}) 
show that the representation 
has the form
\begin{equation}
A=diag\left(
\left(
\begin{array}{cc}
1 & 1 \\
0 & 1\\
\end{array}
\right),q,1\right), \;\;\;B_1=e_{32} \;\;\;B_2=e_{34}.
\label{f75}
\end{equation} 
and therefore $B(A)^2=0$.

3. $\epsilon=q$. By (\ref{f71}) and (\ref{f72}) 
we have $\gamma=0$ and by (\ref{f66}) and (\ref{f67}) 
$\beta=ce_{12}$. With (\ref{f61}) and (\ref{f62})
this implies the representation has the form
\begin{equation}
A=diag\left(
\left(
\begin{array}{cc}
q & 1 \\
0 & q\\
\end{array}
\right),q,1\right), \;\;\;B_1=e_{14} , \;\;\;B_2=e_{34}
\label{f76}
\end{equation}
and again $B(A)^2=0$.\

4. $\epsilon=q^2$. By (\ref{f71}) and (\ref{f72})
 we have $\gamma=0$ and equalities (\ref{f66}) 
and (\ref{f67}) imply $\beta=ce_{11}$. So the representation has the form
\begin{equation}
A=diag\left(
\left(
\begin{array}{cc}
q^2 & 1 \\
0 & q^2\\
\end{array}
\right),q,1\right), \;\;\;B_1=e_{13} , \;\;\;B_2=e_{34}.\label{f77}
\end{equation}
This representation coincides with (\ref{f36}).\hfill$\Box$

Note that if $q^3=1$ or $q^4=1,$ then Theorem \ref{t2} 
is not valid. For $q=e^{\frac{2\pi i}{3}}$ and $q=\pm i$ there exist, 
respectively, three dimensional and four dimensional irreducible 
representations (see for instance \cite{Mi}) while for $q^m\neq 1$ all finite dimensional 
representations of the $q$-spinor are one dimensional (see \cite{Sm}, 
Chapter 2).
\section{Equivalence of representations}
All irreducible finite 
dimensional representations of $GL_q(2,C),$\ $q^{m}\neq 1$ 
are one dimensional.
This is a folklore fact which can be easily obtained from the
Y.S. Soibelman work \cite{So} or from the FRT-duality \cite{RTT}
and the fact that $U_q(g)$ is pointed (see \cite{M93b}).
We need the following corollary from this fact.

\begin{corollary}.
Every finite dimensional representation of $GL_q(2,C),$\ $q^m\neq 1,$ 
is triangular; i.e.  
it is equivalent to a representation by triangular matrices
$a_{ij}\rightarrow A_{ij}.$
\label{c2}
\end{corollary}

\begin{corollary}.
 For every finite dimensional representation 
$a_{ij}\rightarrow A_{ij}$ of $GL_q(2,C), q^m\neq 1$, 
the elements $A_{11}$,
$A_{22}$ are invertible, while $A_{12}$, $A_{21}$ are nilpotent.
\label{c3}
\end{corollary}

{\it Proof}. From Corollary \ref{c2} we can suppose that 
$A_{ij}$ are triangular 
matrices. In this case the matrix 
\begin{equation}
(q-q^{-1})^{-1}(A_{11}A_{22}-A_{22}A_{11})
\end{equation}
has only zero entries on the main diagonal. This matrix is equal to 
$A_{12}A_{21}.$ From this follows that the main diagonal of 
$A_{11}A_{22}$ and that of the invertible matrix $det_q$=
$A_{11}A_{12}-qA_{12}A_{21}$ coincide. 
This means that $A_{11}$ and $A_{22}$
have no zero terms on the main 
diagonal and therefore they are invertible.
\hfill$\Box$\

Now we are ready to answer the question when inner actions defined by two
different representations are equivalent to each other.  
\begin{theorem}.
 Let $a_{ij}\rightarrow A_{ij}$ and 
$a_{ij}\rightarrow A'_{ij}$ be two representations of $GL_q(2,C)$ in
${\it C}(1,3)$. Then Hopf algebra actions
\begin{equation}
a_{ij}\cdot v=\sum_kA_{ik}vA^*_{kj}
\label{f96}
\end{equation}
and
\begin{equation}
a_{ij}*v=\sum_kA'_{ik}vA'^*_{kj}
\label{f97}
\end{equation}
are equivalent if and only if
$$A'_{11}=uA_{11}u^{-1}\alpha_1,\;
A'_{12}=uA_{12}u^{-1}\alpha_2,\;$$
\begin{equation}
A'_{21}=uA_{21}u^{-1}\alpha_1,\;
A'_{22}=uA_{22}u^{-1}\alpha_2,
\label{f98}
\end{equation}
for some nonzero complex numbers $\alpha_1$, $\alpha_2$ and invertible 
$u\in {\it C}(1,3)$. Formulas $(\ref{f96})$ and $(\ref{f97})$ 
define the same action if and only 
if the elements $A'_{ij}$ are connected with $A_{ij}$ 
by formulas $(\ref{f98}),$ with 
$u=1$.
\label{t4}
\end{theorem}

{\it Proof}. First of all, it is easy to see that if 
$a_{ij}\rightarrow A_{ij}$ is a representation and $A'_{ij}$ is defined by
(98), then $a_{ij}\rightarrow A'_{ij}$  is also a representation of
$GL_q(2,C)$. In order to show that they define equivalent actions let us
present formulas(\ref{f96}),(\ref{f97}) in matrix form
\begin{equation}
\left(
\begin{array}{cc}
a_{11}\cdot v & a_{12}\cdot v\\
a_{21}\cdot v & a_{22}\cdot v
\end{array}
\right)=
A\left(
\begin{array}{cc}
v & 0 \\
0 & v
\end{array}
\right)A^{-1}
\label{f99} \end{equation}
and
\begin{equation}
\left(
\begin{array}{cc}
a_{11}*v & a_{12}*v \\
a_{21}*v & a_{22}*v
\end{array}
\right)=
A_1\left(
\begin{array}{cc}
v & 0 \\
0 & v
\end{array}
\right)A^{-1}_1,
\label{f100} \end{equation}
where by Theorem 1
\begin{equation}
A=
\left(
\begin{array}{cc}
A_{11} & A_{12}\\
A_{21} & A_{22}
\end{array}
\right),\;\;\;
A^{-1}=
\left(
\begin{array}{cc}
A^*_{11} & A^*_{12}\\
A^*_{21} & A^*_{22}
\end{array}
\right)
\label{f101} \end{equation}
and
\begin{equation}
A_1=
\left(
\begin{array}{cc}
A'_{11} & A'_{12}\\
A'_{21} & A'_{22}
\end{array}
\right)=
\left(
\begin{array}{cc}
uA_{11}u^{-1}\alpha_1 & uA_{12}u^{-1}\alpha_2\\
uA_{21}u^{-1}\alpha_1 & Au_{22}u^{-1}\alpha_2
\end{array}
\right)=
\label{f102} \end{equation}
$$
\left(
\begin{array}{cc}
u & 0\\
0 & u
\end{array}
\right)A
\left(
\begin{array}{cc}
u^{-1} & 0\\
0 & u^{-1}
\end{array}
\right)
\left(
\begin{array}{cc}
\alpha_1 & 0\\
0     & \alpha_2
\end{array}
\right).
$$
From this formula we have
\begin{equation}
A^{-1}_1=
\left(
\begin{array}{cc}
\alpha^{-1}_1 & 0\\
0 & \alpha^{-1}_2
\end{array}
\right)
\left(
\begin{array}{cc}
u & 0\\
0 & u
\end{array}
\right)A^{-1}
\left(
\begin{array}{cc}
u^{-1} & 0\\
0 & u^{-1}
\end{array}
\right).
\label{f103} \end{equation}
Taking into account that $diag(v,v)$ commute with all $2\times 2$ matrix with
complex coefficients (and in particular with $diag(\alpha_1,\alpha_2)$) we
have
\begin{equation}
A_1
\left(
\begin{array}{cc}
v  &  0 \\
0  &  v
\end{array}
\right)A^{-1}_1=
\label{f104} \end{equation}
$$
\left(
\begin{array}{cc}
u  &  0 \\
0  &  u
\end{array}
\right)A
\left(
\begin{array}{cc}
u^{-1}  &  0 \\
0  &  u^{-1}
\end{array}
\right)
\left(
\begin{array}{cc}
v  &  0 \\
0  &  v
\end{array}
\right)
\left(
\begin{array}{cc}
u  &  0 \\
0  &  u
\end{array}
\right)A^{-1}
\left(
\begin{array}{cc}
u^{-1}  &  0 \\
0  &  u^{-1}
\end{array}
\right)=
$$
$$
\left(
\begin{array}{cc}
u  &  0 \\
0  &  u
\end{array}
\right)A
\left(
\begin{array}{cc}
u^{-1}vu  &  0 \\
0  &  u^{-1}vu
\end{array}
\right)A^{-1}
\left(
\begin{array}{cc}
u^{-1}  &  0 \\
0  &  u^{-1}
\end{array}
\right)=
$$
$$
\left(
\begin{array}{cc}
ua_{11}\cdot(u^{-1}vu)u^{-1}   & ua_{12}\cdot(u^{-1}vu)u^{-1}\\
ua_{21}\cdot(u^{-1}vu)u^{-1}   & ua_{22}\cdot(u^{-1}vu)u^{-1}
\end{array}
\right).
$$
If we set $w=u^{-1}vu,$ then by(\ref{f100}) and(\ref{f104}) we will get
\begin{equation}
a_{ij}*(uwu^{-1})=u(a_{ij}\cdot w)u^{-1}.
\label{f105} \end{equation}
The last formula concides with(\ref{fequiv}) and therefore the actions are equivalent.
If $u=1,$ then(\ref{f105}) shows that the actions coincide.\

Inversely, let $a_{ij}\rightarrow A_{ij}$, $a_{ij}\rightarrow A'_{ij}$ be two
representations which define equivalent (or equal) actions. Then we have the
relations of type(\ref{fequiv}) (with $u=1,$ if the actions coincide). In the matrix from
for $v=uwu^{-1}$ this relations are equivalent to
\begin{equation}
A_1 
\left(
\begin{array}{cc}
v  &  0 \\
0  &  v
\end{array}
\right)A^{-1}_1=
\label{f106} \end{equation}
$$
\left(
\begin{array}{cc}
u  &  0 \\
0  &  u
\end{array}
\right)A
\left(
\begin{array}{cc}
u^{-1}vu  &  0 \\
0  &  u^{-1}vu
\end{array}
\right)A^{-1}
\left(
\begin{array}{cc}
u^{-1}  &  0 \\
0  &  u^{-1}
\end{array}
\right)
$$
Let us multiply this relation from the left by 
$$\left(
\begin{array}{cc}
u  &  0 \\
0  &  u
\end{array}
\right)A^{-1}
\left(
\begin{array}{cc}
u^{-1}  &  0 \\
0  &  u^{-1}
\end{array}
\right)$$
and from the right by $A_1$. We will get
\begin{equation}
\left(
\begin{array}{cc}
u  &  0 \\
0  &  u
\end{array}
\right)A^{-1}
\left(
\begin{array}{cc}
u^{-1}  &  0 \\
0  &  u^{-1}
\end{array}
\right)A_1
\left(
\begin{array}{cc}
v  &  0 \\
0  &  v
\end{array}
\right)=
\label{f107} \end{equation}
$$
\left(
\begin{array}{cc}
v  &  0 \\
0  &  v
\end{array}
\right)
\left(
\begin{array}{cc}
u  &  0 \\
0  &  u
\end{array}
\right)A^{-1}
\left(
\begin{array}{cc}
u^{-1}  &  0 \\
0  &  u^{-1}
\end{array}
\right)A_1.
$$
If we denote by $\alpha_{ij}$ the coefficients of
$\left(
\begin{array}{cc}
u  &  0 \\
0  &  u
\end{array}
\right)A^{-1}
\left(
\begin{array}{cc}
u^{-1}  &  0 \\
0  &  u^{-1}
\end{array}
\right)A_1$, then this relation reduces to
\begin{equation}
\left(
\begin{array}{cc}
\alpha_{11}v  &  \alpha_{12}v \\
\alpha_{21}v  &  \alpha_{22}v
\end{array}
\right)=
\left(
\begin{array}{cc}
v\alpha_{11}  &  v\alpha_{12} \\
v\alpha_{21}  &  v\alpha_{22}
\end{array}
\right),
\label{f108} \end{equation}
i.e. $[\alpha_{ij},v]=0$. As $v$ is an arbitrary element from ${\it C}(1,3)$, the
elements $\alpha_{ij}$ belong to the center of ${\it C}(1,3)$. The center of
${\it C}(1,3)$ coincide with the set of all complex numbers , i.e. 
$\alpha_{ij}\in C$. Thus, we have
\begin{equation}
A_1=\left(
\begin{array}{cc}
u^{-1}  &  0 \\
0  &  u^{-1}
\end{array}
\right)A
\left(
\begin{array}{cc}
u  &  0 \\
0  &  u
\end{array}
\right)
\left(
\begin{array}{cc}
\alpha_{11}  &  \alpha_{12} \\
\alpha_{21}  &  \alpha_{22}
\end{array}
\right)
\label{f109} \end{equation}
or in details
\begin{equation}
A'_{11}=u^{-1}(A_{11}\alpha_{11}+A_{12}\alpha_{21})u,\;\;
A'_{22}=u^{-1}(A_{21}\alpha_{12}+A_{22}\alpha_{22})u
\label{f110} \end{equation}
\begin{equation}
A'_{12}=u^{-1}(A_{11}\alpha_{12}+A_{12}\alpha_{22})u,\;\;
A'_{21}=u^{-1}(A_{21}\alpha_{11}+A_{22}\alpha_{21})u.
\label{f111} \end{equation}
We know that $a_{ij}\rightarrow A'_{ij}$ is a representation of $GL_q(2,C)$
 and therefore by Corollary 3 the elements $A'_{12}$, $A'_{21}$ are
nilpotent. By Corollary 2 we can suppose that 
$A_{ij}$ 
are triangular matrices of the form (90) and by Corollary 3 all entries on the 
main diagonal of $A_{11}$, $A_{22}$ are nonzero, while all entires on the main
diagonal of $A_{12}$, $A_{21}$ are zero. Now, if $\alpha_{12}\neq 0$ then the
first formula of(\ref{f111}) shows that all entries on the main diagonal of 
$uA'_{12}u^{-1}$ are nonzero, therefore $uA'_{12}u^{-1}$ as well as $A'_{12}$
are invertible. This is impossible as
$A'_{12}$ is nilpotent. So $\alpha_{12}=0$.\

In the same way $\alpha_{21}=0$. Thus, if we denote $\alpha_1=\alpha_{11}$,
$\alpha_2=\alpha_{22}$ then(\ref{f110}) and(\ref{f111}) 
reduce to(\ref{f98}).\hfill$\Box$

{\bf Definition}. {\it We call two representations $a_{ij}\rightarrow A_{ij}$, 
$a_{ij}\rightarrow A'_{ij}$ equivalent if they define equivalent actions of 
$GL_q(2,C)$; i.e., if relations $(\ref{f98})$ are satisfied}.\\

As a 
homomorphism $\varphi: GL_q(2,C)\rightarrow {\it C}(1,3)$ defines a four 
dimensional left module over the algebra $GL_q(2,C)$ and viceversa,
Theorem \ref{t4} shows that the equivalence of representations
means that corresponding modules $V_1$, $V_2$ are related by formula 
$V_1\simeq V_2\otimes {\it U}$, where ${\it U}$ is any one dimensional
module.\

Indeed, numbers $\alpha_1$, $\alpha_2$ define one dimensional representation
$a_{11}\rightarrow \alpha_1$, $a_{12}\rightarrow 0$, $a_{21}\rightarrow 0$,
$a_{22}\rightarrow \alpha_2$ and every one dimensional representation
$a_{ij}\rightarrow \alpha_{ij}$ is defined by two nonzero numbers
$\alpha_{11}$, $\alpha_{22}$, while $\alpha_{12}$=$\alpha_{21}$=$0$
($\alpha_{11}\alpha_{12}$=$q\alpha_{12}\alpha_{11}$=$q\alpha_{11}\alpha_{12}$,
$\alpha_{11}\neq 0$). 

\section{Algebra $R_q$ and its representations}
Let us consider some simplification 
of the algebraic structure of $GL_q(2,C).$
Let the algebra $R_q$ be generated by the
elements $a_{11}, a_{12}, a_{21}, r_{22}$
with the relations between $a_{11}, a_{12}, a_{21}$ defined 
by (\ref{p1}) and
the additional relations
\begin{equation}
a_{12}r_{22}=qr_{22}a_{12}, \ \ \ a_{21}r_{22}=qr_{22}a_{21},
\label{sp1}
\end{equation}
\begin{equation}
a_{11}r_{22}=r_{22}a_{11}.
\label{sp2}
\end{equation}
In fact this means that $R_q$ is defined by the same system of spinors with
commutative diagonal:
\begin{equation} 
\begin{picture}(80,80)
\put(15.29,2.39){\makebox(0,0)[cc]{$a_{21}$}}
\put(60.0,2.59){\makebox(0,0)[cc]{$r_{22}$}}
\put(60.0,45.46){\makebox(0,0)[cc]{$a_{12}$}}
\put(15.29,45.46){\makebox(0,0)[cc]{$a_{11}$}}
\put(22,45.46){\vector(1,0){28}} 
\put(22,2.39){\vector(1,0){28}}
\put(15.29,38){\vector(0,-1){28}}
\put(58.32,38){\vector(0,-1){28}}
\put(20.0,8.0){\line(1,1){30}}
\put(25.0,38.0){\line(1,-1){30}}
\end{picture}
\label{p1}
\end{equation}
\begin{theorem}.
If $a_{ij}\rightarrow A_{ij}$ is a finite dimensional
representation of the algebra $GL_q(2,C)$, $q^m\neq 1$, then
\begin{equation}
a_{ij}\rightarrow A_{ij}, \;\;
r_{22}\rightarrow R_{22}=A_{22}-A_{12}A^{-1}_{11}A_{21}
\label{repspi}
\end{equation}
is a representation of $R_q$ with invertible $A_{11}R_{22}$. Inversely, if
$a_{ij}\rightarrow A_{ij}$ $(i\neq2$ or  $j\neq 2)$, $r_{22}\rightarrow R_{22}$
is a finite dimensional representation of $R_q$, such that $A_{11}$, $R_{22}$
are invertible then
\begin{equation}
a_{ij}\rightarrow A_{ij},\;\;
a_{22}\rightarrow A_{22}=R_{22}+A_{12}A^{-1}_{11}A_{21}
\label{repglq}
\end{equation}
is a representation of algebra $GL_q(2,C)$.
 In this case $(\ref{repglq})$ is a representation of $SL_q(2,C)$ iff 
$R_{22}=A^{-1}_{11}$.
\label{spin}
\end{theorem}

{\it Proof}. Let $a_{ij}\rightarrow A_{ij}$ be a finite dimensional 
representation of $GL_q(2,C)$. Then by Corollary \ref{c3}
the elements $A_{11}$, $A_{22}$ are invertible, therefore formula 
(\ref{repspi}) is correctly defined. 
It is easy to see that
\begin{equation}
[A_{12}, R_{22}]_q=[A_{21}, R_{22}]_q=[A_{11},R_{22}]=0,
\label{sp3}
\end{equation}
where we have denoted $[x,y]_q=xy-qyx, \ [x,y]=xy-yx.$

Besides we have
$$
R_{22}=A_{22}-qA^{-1}_{11}A_{12}A_{21}=
A^{-1}_{11}(A_{11}A_{22}-qA_{12}A_{21})=A^{-1}_{11}det_q,
$$
which shows that $R_{22}$ is invertible.\

Inversely, let $r_{22}\rightarrow R_{22}$, $a_{ij}\rightarrow A_{ij}$; 
$i\neq 2$ or $j\neq 2$, be a finite dimensional representation of $R_q$ with 
invertible $A_{11}$, $R_{22}$. Straightforward calculations show that
\begin{equation}
[A_{12}, A_{22}]_q=[A_{21}, A_{22}]_q=0,\;\;
[A_{11}, A_{22}]=(q-q^{-1})A_{12}A_{21}.
\label{fo27}
\end{equation}
 
We have 
$$
det_q=A_{11}A_{22}-qA_{12}A_{21}=A_{11}(R_{22}+A_{12}A^{-1}_{11}A_{21})-
qA_{12}A_{21}=
$$
$$
A_{11}R_{22}+qA_{12}A_{21}-
qA_{12}A_{21}=A_{11}R_{22},
$$
therefore $det_q$ is invertible. \hfill$\Box$

\section{Actions defined by $SL_q(2,C)$-representations}  
For a classification of the  algebra representations of
$GL_q(2,C)$ in ${\it C}(1,3),$ by Theorem 1, it is enough to 
describe representations of
$R_q$ with invertible $a_{11}$, $r_{22}$. 
In this case we have two $q$-spinors
($a_{11}, a_{12}$), $(a_{11}, a_{21})$ and two $q^{-1}$-spinors 
$(r_{22}, a_{12})$, $(a_{22}, a_{21})$ with invertible first terms, such that
$a_{11}r_{22}=r_{22}a_{11}$.\

Note that if $(x,y)$ is a $q$-spinor then $(x^{-1},y)$ is a $q^{-1}$-spinor.
Therefore for each of the two representations $(x, y)\rightarrow (A_{11}, A_{12}),$
\ $(x, y)\rightarrow (A_{11}, A_{21})$ with invertible
$A_{11}$ and commuting
$A_{12}$, $A_{21}$ we have a representation of the algebra $R_q:$
\begin{equation}
a_{11}\rightarrow A_{11},\;\;
a_{12}\rightarrow A_{12},\;\;
a_{21}\rightarrow A_{21},\;\;
r_{22}\rightarrow A^{-1}_{11}.
\label{repslq}
\end{equation}
By Theorem \ref{spin} this representation defines
a representation of $SL_q(2,C)$.
In this way, using $q$-spinor representations given in
 Theorem \ref{t2} ,
we can write a number of $SL_q(2,C)$-representations  with
nontrivial ``perturbation" --- see all representations in the Table,
which are marked by the letter $S$: $S1, S2a, S2a^{\prime },$ e.t.c.
\begin{theorem}.
Every representation $a_{ij}\rightarrow A_{ij}$ of the algebra 
$SL_q(2,C), \ q^m\neq 1$ in $C(1,3)$ with nontrivial ``perturbation"
is equivalent to one of the representations marked by letter $S$ 
in the Table. Invariants and operator algebras of corresponding inner
actions are given in the Table.
\label{slq}
\end{theorem}
 
{\it Proof}. Let $a_{ij}\rightarrow A_{ij}$ be a representation 
of $ GL_q(2,C)$ and 
$a_{ij}\rightarrow A_{ij}$, $r_{22}\rightarrow R_{22}$ be the corresponding 
representation of $R_q$; i.e. $R_{22}$=$A_{22}-qA^{-1}_{11}A_{12}A_{21}$ 
(see (\ref{repspi})). By Theorem 1 we have $R_{22}$=$A^{-1}_{11}$ and therefore
\begin{equation}
A_{22}=A^{-1}_{11}(1+qA_{12}A_{21}).
\label{s}
\end{equation}
We know that $a_{11}\rightarrow A_{11}$, $a_{21}\rightarrow A_{21}$ and
$a_{11}\rightarrow A_{11}$, $a_{12}\rightarrow A_{12}$ are two representations
of the $q$-spinor with an invertible first term. Recall that for a matrix 
$A$ we have denoted by $B(A)$ the linear space of all 
matrices $B$ such that
$AB=qBA$. Thus 
$A_{12}, A_{21}\in B(A_{11})$ and so $B(A_{11})^2\neq 0$. 
By Theorem \ref{t2}  we 
can assume that $A_{11}$ is one of the seven matrices 
which appear in this theorem. 
Let us consider the seven cases separately.\

{\bf 1. $\bf A_{11}=diag(q^2, q, 1, 1)$.} In this case we have
\begin{equation} 
A_{12}=\alpha e_{12}+\beta e_{23}+\gamma e_{24},\;\;\; 
A_{21}=\alpha_1 e_{12}+\beta_1 e_{23}+\gamma_1 e_{24},
\label{s1.1}
\end{equation}
therefore
\begin{equation}
A_{12}A_{21}=\alpha\beta_1e_{13}+\alpha\gamma_1e_{14}=A_{21}A_{12}=
\alpha_1\beta e_{13}+\alpha_1\gamma e_{14}\neq 0.
\label{s1.2}
\end{equation}
These relations imply that
\begin{equation}
0=\alpha\beta_1-\alpha_1\beta=det
\left(
\begin{array}{cc}
\alpha   & \beta \\
\alpha_1 & \beta_1
\end{array}
\right),\;\;
0=\alpha\gamma_1-\alpha_1\gamma=det
\left(
\begin{array}{cc}
\alpha   & \gamma \\
\alpha_1 & \gamma_1
\end{array}
\right).
\label{s1.3}
\end{equation}
The first of equalities (\ref{s1.3}) shows that $(\alpha, \beta)$ and 
$(\alpha_1, \beta_1)$ are linearly dependent vectors and by (\ref{s1.2})
we can write
$(\alpha_1, \beta_1)$=$\frac{\alpha_1}{\alpha}(\alpha,\beta)$. In the same
way $(\alpha_1, \gamma_1)$=$\frac{\alpha_1}{\alpha}(\alpha, \gamma)$. These
two relations are equivalent to
\begin{equation}
A_{21}=\epsilon A_{12},\;\;\epsilon=\frac{\alpha_1}{\alpha}\neq 0.
\label{s1.4}
\end{equation}
Let us consider a matrix of the form
$
{\it U}=diag(1, d, M),
$
where $M$ is an invertible $2\times 2$ matrix and $d$ a nonzero complex 
number. This matrix commutes with $A_{11}$ therefore the conjugation by
it will not change $A_{11}$, while $A_{12}$ is changed in the following
way
\begin{equation}
{\it U}A_{12}{\it U}^{-1}=
\left(
\begin{array}{ccc}
0  & \alpha d^{-1} &   \begin{array}{cc}0\ \   &\ \   0\end{array}  \\
0  &    0   &  (d\beta,\alpha\gamma)  M^{-1} \\
0  &    0   & \begin{array}{cc}0\ \   &\ \   0\end{array} \\
0  &    0   & \begin{array}{cc}0\ \   &\ \   0\end{array}
\end{array}
\right).
\label{s1.7}
\end{equation}
Let us take $d=\alpha$. Then $(d\beta, \alpha \gamma )$ =
$(\alpha\beta, \alpha\gamma)\neq 0$ and
there exists an invertible $2\times 2$ matrix $M$ such that
$(\alpha\beta, \alpha\gamma)M^{-1}=(1,0)$. If we replace $\epsilon$ by
$\alpha$, then we will get the representation $S1$ of the Table.\

Let us calculate $\Re $ and the invariants of the $SL_q(2,C)$-action defined by
this representation.\

All three degrees of $A_{11}$ are linearly independent, as the Vandermond
determinant
does not vanish. This means that the algebra $\Re $ contains the elements 
$e_{11}, e_{22}, e_{33}+e_{44}$. We have $A_{12}$=$e_{12}+e_{23}\in \Re $
and $e_{11}A_{12}$=$e_{12}\in \Re $. Thus $e_{12}$, $e_{23}$, $e_{13}$=
$e_{12}e_{23}\in \Re $. The linear space generated by these six elements is 
a subalgebra given in S1 in the Table, which evidentely
is  isomorphic to the algebra of triangular
$3 \times 3$ matrices $T_3.$ 
For calculating the invariants we can use Lemma 1.\
It is easy to see that the centralizer of $\Re $ 
is equal to the algebra of matrices 
$\beta E+\gamma e_{44}+\delta e_{43},$
which is isomorphic to the algebra of $2 \times 2$ 
triangular matrices.
In fact in this case we have just two basic nontrivial invariants 
(see  \ref{matrix})
\begin{equation}
I_1=(1-\gamma_0)(-\gamma_1+i\gamma_2)\gamma_3,\;\;
I_2=(1-\gamma_0)(1-i\gamma_{12}),
\label{s1.12}
\end{equation}
while all others are linear combinations of them and of the unit element.\

{\bf 2. $A_{11}=diag(q^2,q,q,1).$} By theorem \ref{t2} we have $A_{12}$=
$\alpha e_{12}+\beta e_{13}+\gamma e_{24}+\delta e_{34}$, $A_{21}$=
$\alpha_1 e_{12}+\beta_1 e_{13}+\gamma_1 e_{24}+\delta_1 e_{34}$.\
Therefore
\begin{equation}
A_{12}A_{21}=\left(\alpha\gamma_1+\beta\delta_1\right)e_{14}=A_{21}A_{12}=
\left(\alpha_1\gamma+\beta_1\delta\right)e_{14}\neq 0.
\label{s2.1}
\end{equation}

From (\ref{s2.1}) we have
\begin{equation}
\alpha\gamma_1+\beta\delta_1=\alpha_1\gamma+\beta_1\delta\neq 0.
\label{s2.2}
\end{equation}
The matrix $A_{11}$ commutes with all matrices ${\it U}$=$diag(d, M, 1)$, 
where $d$ is a nonzero complex number and $M$ is an invertible $2\times 2$ 
matrix. Conjugation by a matrix ${\it U}$ does not change $A_{11}$, while 
$A_{12}$, $A_{21}$ are changed by the following formulae
$$A_{12}\rightarrow
\left(
\begin{array}{ccc}
0 & (d\alpha , d\beta )M^{-1} & 0 \\
\begin{array}{c} 0 \\ 0\end{array} & 
\begin{array}{cc} 0 \ \ &\ \  0 \\ 0\ \  & \ \ 0\end{array} &
M\left(\begin{array}{c} \gamma \\ \delta \end{array}\right) \\   
0 &  0\ \ \ \ \, \ \ 0 & 0   
\end{array}
\right);$$
\begin{equation}
\ \ \ A_{21}\rightarrow
\left(
\begin{array}{ccc}
0 & (d\alpha_1 , d\beta_1 )M^{-1} & 0 \\
\begin{array}{c} 0 \\ 0\end{array} & 
\begin{array}{cc} 0 \ \ &\ \  0 \\ 0\ \  &\ \  0\end{array} &
M\left(\begin{array}{c} \gamma_1 \\ \delta_1\end{array}\right) \\   
0 &  0\ \ \, \ \ \ \ 0 & 0   
\end{array}
\right).
\label{Suemi}
\end{equation}
In particular, the matrix 
$\left(
\begin{array}{cc}
\gamma & \gamma_1 \\
\delta & \delta_1
\end{array}
\right)$
under this conjugation is multiplied by $M$ from the left. 
For this matrix we have two possibilities.\

a) $det\left(\begin{array}{cc}
\gamma & \gamma_1 \\
\delta & \delta_1
\end{array}
\right)\neq 0$. In this case this matrix is invertible and we can take $M$ to
be its inverse
and (\ref{Suemi}) reduces to the form
$A_{12}\rightarrow d\alpha e_{12}+d\beta e_{13}+e_{24},$\
$A_{21}\rightarrow d\alpha _1e_{12}+d\beta _1e_{13}+e_{34},$\
with new parameters $\alpha$, $\beta$, $\alpha_1$, $\beta_1$. Relations 
(\ref{s2.2})
become  $\alpha\cdot 0+\beta\cdot 1\neq 0$; i.e. $\beta\neq 0$, and
$\beta$=$\alpha_1\cdot 1+\beta_1\cdot 0$=$\alpha_1\neq 0$. If we let
$d$=$\beta^{-1}$=$\alpha^{-1}_1,$ then changing 
$\alpha^{-1}_1\alpha\rightarrow\alpha$, 
$\alpha^{-1}_1\beta_1\rightarrow\beta ,$ we get the representation $S2a.$\

For calculating the algebra $\Re $ we can make the same procedure as in the
first case. The first three degrees of $A_{11}$ are linearly independent and
so $\Re $ contains $e_{11}$, $e_{22}+e_{33}$, $e_{44}$. Multiplication of
$A_{12}$ and $A_{21}$ by $e_{44}$ from the right gives two elements $e_{24}$,
$e_{34}$ and also $\alpha e_{12}+e_{13}$, $e_{12}+\beta e_{13}$. 

If 
$\alpha\beta\neq 1$ then these two elements are linearly independent and
$e_{12}$, $e_{13}\in \Re $. So $\Re $ consists of matrices of the form 
given in the Table for representation $S2a$.\

If $\alpha\beta=1$ then the elements $\alpha e_{12}+e_{13}$ and 
$e_{12}+\beta e_{13}$ are linearly dependent and all matrices from $\Re $ 
have $(-_{12}, -_{13})$-entries proportional to 
$(\alpha , 1)=\alpha (1, \beta )$; i.e. $\Re $ has
the form given in the Table for $S2a^{\prime }$. 
Easy calculations show that in both cases the centralizer of 
$\Re $ is equal to $1\cdot C$.\

b) $det\left(
\begin{array}{cc}
\gamma & \gamma_1 \\
\delta & \delta_1
\end{array}
\right)=0$. By the relation (\ref{s2.1}) we have 
$\left(
\begin{array}{c}
\gamma\\
\delta
\end{array}
\right)\neq 0$,
$\left(
\begin{array}{c}
\gamma_1\\
\delta_1
\end{array}
\right)\neq 0,$ therefore it is possible to find a matrix $M_1$ such that
\begin{equation}
M_1
\left(
\begin{array}{c}
\gamma \\
\delta
\end{array}
\right)=
\left(
\begin{array}{c}
1\\
0
\end{array}
\right);\ \ 
M_1
\left(
\begin{array}{c}
\gamma_1 \\
\delta_1
\end{array}
\right)=
\left(
\begin{array}{c}
\gamma'_1\\
0
\end{array}
\right),\ \ 
\gamma'_1\neq 0.
\label{s2.6}
\end{equation}
Let $M=M_2M_1$, where $M_2$ is an invertible matrix, such that 
$M_2\left(
\begin{array}{c}
1 \\
0
\end{array}
\right)$=
$\left(
\begin{array}{c}
1 \\
0
\end{array}
\right)$; i.e.
$
M_2=
\left(
\begin{array}{cc}
1 & u' \\
0 & v
\end{array}
\right),\ \ 
v\neq 0.
\label{s2.7}
$
In this case formula (\ref{Suemi}) takes the form
$$A_{12}\rightarrow
\left(
\begin{array}{ccc}
0 & (d\alpha' , d\beta' )M^{-1}_2 & 0 \\
\begin{array}{c} 0 \\ 0\end{array} & 
\begin{array}{cc} 0 \ \ &\ \  0 \\ 0\ \  & \ \ 0\end{array} &
\begin{array}{c} 1 \\ 0 \end{array} \\   
0 &  0\ \ \ \ \, \ \ 0 & 0   
\end{array}
\right);$$
\begin{equation}
\ \ \ A_{21}\rightarrow
\left(
\begin{array}{ccc}
0 & (d\alpha'_1 , d\beta'_1 )M^{-1}_2 & 0 \\
\begin{array}{c} 0 \\ 0\end{array} & 
\begin{array}{cc} 0 \ \ &\ \  0 \\ 0\ \  &\ \  0\end{array} &
\begin{array}{c} \gamma'_1 \\ 0\end{array} \\   
0 &  0\ \ \, \ \ \ \ 0 & 0   
\end{array}
\right),
\label{sue}
\end{equation}
where $(\alpha', \beta')$=$(\alpha, \beta)M^{-1}_1$,
$(\alpha'_1, \beta'_1)$=$(\alpha_1, \beta_1)M^{-1}_1$ and
\begin{equation}
M^{-1}_2=
\left(
\begin{array}{cc}
1 & u \\
0 & v^{-1}
\end{array}
\right),\;\; u=-u'v^{-1}.
\label{s2.8}
\end{equation}
Shortly we can write (\ref{sue}) in the form
\begin{equation}
A_{12}\rightarrow d\alpha''e_{12}+d\beta''e_{13}+e_{24},\;\;
A_{21}\rightarrow d\alpha''_1e_{12}+d\beta''_1e_{13}+\gamma''_1e_{24},
\label{suem}
\end{equation}
where
\begin{equation}
(\alpha'', \beta'')=(\alpha', \beta')
\left(
\begin{array}{cc}
1 & u \\
0 & v^{-1}
\end{array}
\right)=
(\alpha', \alpha'u+\beta'v^{-1})
\label{s2.9}
\end{equation}
and
\begin{equation}
(\alpha''_1, \beta''_1)=(\alpha'_1, \beta'_1)
\left(
\begin{array}{cc}
1 & u \\
0 & v^{-1}
\end{array}
\right)=
(\alpha'_1, \alpha'_1u+\beta'_1v^{-1}).
\label{s2.10}
\end{equation}
By using (\ref{s2.2}), applied to primed parameters, we have
\begin{equation}
\alpha'\gamma_1'+\beta'\cdot 0=\alpha'_1\cdot 1+\beta'_1\cdot 0\neq 0;
\mbox{ i.e.}\;\; \alpha'\neq 0,
\;\;\alpha'\gamma'_1=\alpha'_1.
\label{s2.11}
\end{equation}
If we suppose $u=-(\alpha')^{-1}\beta'v^{-1}$ then we get $\beta''=0$
and 
\begin{equation}
\beta''_1=-\alpha'_1(\alpha')^{-1}\beta'v^{-1}+\beta'_1v^{-1}=
(\beta'_1-\gamma'_1\beta')v^{-1}. \label{last}
\end{equation}
If $\beta'_1\neq \gamma'_1\beta'$ in formula (\ref{last}), 
then we can take $v$=
$(\beta'_1-\gamma'_1\beta')(\alpha')^{-1}$. In this case $\beta''_1$=$\alpha'$ and
we can let $d$=$(\alpha')^{-1}$ in (\ref{suem}) 
in order to obtain $S2b$
(where by $\alpha$ we mean $(\alpha')^{-1}\alpha''_1$=
$(\alpha')^{-1}\alpha'_1=\gamma'_1)$.\

If $\beta'_1$=$\gamma_1\beta'$ in (\ref{last}), then $\beta''_1$=$0$ and (\ref{suem}) with 
$d$=$(\alpha')^{-1}$ is equal to $S2b^{\prime }.$ In this case the 
algebra $\Re $ contains the elements $e_{11}$, $e_{22}+e_{33}$, $e_{44}$, $e_{12}$,
$e_{24}$, $e_{14}$ which form a basis of this algebra. This 
algebra is isomorphic to $T_3$ 
(see in the Table an action S2$b^{\prime }$.)

The centralizer of $\Re $ is equal to the set of diagonal 
matrices of the form
$diag(\beta, \beta,\delta,\beta);$ therefore it is isomorphic to the 
direct sum $C\oplus C.$ 
It is interesting to note that in this case the 
algebra $\Re $ is abstractly 
isomorphic to the algebra $\Re $ for representation $S1$, 
but they are not conjugate in ${\it C}(1,3)$ because 
they have nonisomorphic centralizers.\

Thus, we have that the action of the 
quantum group $SL_q(2,C)$ defined
by the representation $S2b^{\prime }$ has only one basic invariant; i.e.
\begin{equation}
I=(1-\gamma_0)(1+i\gamma_{12}).
\end{equation}

For the representation $S2b$ the algebra $\Re $ 
contains one more matrix, $e_{13}$,
 therefore it is a 7-dimensional algebra (see the Table), 
whose centralizer is equal to $C\cdot 1$ and 
the corresponding $SL_q(2,C)$-action has only trivial
invariants.

{\bf 3. $\bf A_{11}=diag(q^2, q^2, q, 1)$.} By theorem \ref{t2} we have
\begin{equation}
A_{12}=\alpha e_{13}+\beta e_{23}+\gamma e_{34},\;\;
A_{21}=\alpha_1e_{13}+\beta_1e_{23}+\gamma_1e_{34}\; ;
\end{equation}
thus
\begin{equation}
A_{12}A_{21}=\alpha\gamma_1e_{14}+\beta\gamma_1e_{24}=A_{21}A_{12}=
\alpha_1\gamma e_{14}+\beta_1\gamma e_{24}\neq 0.
\label{s3.1}
\end{equation}
These relations imply
\begin{equation}
0=\alpha\gamma_1-\alpha_1\gamma=det
\left(
\begin{array}{cc}
\alpha & \gamma \\
\alpha_1 & \gamma_1
\end{array}
\right),\;\;
0=\beta\gamma_1-\gamma\beta_1=det
\left(
\begin{array}{cc}
\beta & \gamma \\
\beta_1 & \gamma_1
\end{array}
\right),
\label{s3.2}
\end{equation}
\begin{equation}
(\alpha , \beta )\neq 0,\;\;
\gamma_1\neq 0,\;\;
(\alpha_1 , \beta_1)\neq 0,\;\;
\gamma\neq 0.
\label{s3.3}
\end{equation}
From these relations we have
\begin{equation}
A_{21}=\epsilon A_{12},\;\;
\epsilon=\gamma_1/\gamma\neq 0.
\label{s3.4}
\end{equation}
Let us consider a matrix of the form
$
{\it U}=diag(M, d, 1),
$
where $M$ is an invertible $2\times 2$ matrix and $d$ is a nonzero complex
number. The conjugation by this matrix does not change $A_{11}$ but $A_{12}$
changes in the following way
\begin{equation}
{\it U}A_{12}{\it U}^{-1}=\left(
\begin{array}{cccc}
\begin{array}{c} 0 \\ 0 \end{array} & \begin{array}{c} 0 \\ 0 \end{array} &
M\left( \begin{array}{c} \alpha d^{-1}\\ \beta d^{-1} \end{array}\right)& 
\begin{array}{c} 0 \\ 0 \end{array} \\
0    &    0   &     0      &   d\gamma  \\
0    &    0   &     0      &   0
\end{array}
\right).
\label{3.suemi}
\end{equation}
If we take $d=\gamma^{-1}$ then 
$(\alpha d^{-1}, \beta d^{-1})=(\alpha\gamma , \beta\gamma )$ 
is a nonzero vector. Thus there exists an invertible matrix $M$
such that 
$M\left(
\begin{array}{c}
\alpha \gamma \\
\beta \gamma
\end{array}
\right)$=$
\left(
\begin{array}{c}
1 \\
0
\end{array}
\right)$. In this way we have obtained the representation $S3,$ 
where by $\alpha$ is denoted the parameter $\epsilon ,$ see (\ref{s3.4}).\

In this case the algebra $\Re $ is generated by the elements $e_{11}+e_{22}$,
$e_{33}$, $e_{44}$, $e_{13}$, $e_{34}$, $e_{14}$; i.e. 
it consists of matrices
\begin{equation}
\left(
\begin{array}{cccc}
\epsilon & 0 & * & * \\
0    & \epsilon & 0 & 0 \\
0   &    0   & *  & * \\
0   &   0  &   0  & *
\end{array}
\right).
\end{equation}
This algebra is isomorphic to $T_3,$ 
while its centralizer is the algebra of matrices of the form
$\beta E+\gamma e_{22}+\delta e_{12}$
this is isomorphic to the algebra of triangular $2\times 2$ matrices
$T_2.$ 
Now formulae (\ref{matrix}) show that the 
action defined by this representation has the
following basic invariants
\begin{equation}
I_1=(1+\gamma_1)(1-i\gamma_{12}),
I_2=(1+\gamma_0)(\gamma_1+i\gamma_2)\gamma_3.
\end{equation}
{\bf 4. $\bf A_{11}=diag(q^3, q^2, q, 1)$.} By theorem \ref{t2} we have
\begin{equation}
A_{12}=\alpha e_{12}+\beta e_{23}+\gamma e_{34}, \ \  
A_{21}=\alpha_1 e_{12}+\beta_1 e_{23}+\gamma_1 e_{34}\; ;
\end{equation}
thus
\begin{equation}
A_{12}A_{21}=
\alpha\beta_1e_{13}+\beta\gamma_1e_{24}=A_{21}A_{12}=\alpha_1\beta
e_{13}+\beta_1\gamma e_{24}\neq 0.
\label{s4.1}
\end{equation}
These relations imply
\begin{equation}
0=\alpha\beta_1-\alpha_1\beta=
det\left(
\begin{array}{cc}
\alpha & \beta \\
\alpha_1 & \beta_1
\end{array}
\right),\;\;
0=\beta\gamma_1-\beta_1\gamma=
det\left(
\begin{array}{cc}
\beta & \gamma \\
\beta_1 & \gamma_1
\end{array}
\right). 
\label{s4.2}
\end{equation}
Therefore the vectors $(\alpha , \beta )$ and
$(\alpha_1 , \beta_1 )$ are linearly dependent and so are 
$(\beta , \gamma )$ and
$(\beta_1 , \gamma_1 ).$ By (\ref{s4.1}), one of the numbers $\beta , \beta_1$ 
is nonzero. If, for example, $\beta\neq 0$ then
$(\alpha_1 , \beta_1 )=\frac{\beta_1}{\beta}
(\alpha , \beta )$ and $(\beta_1 , \gamma_1 )=\frac{\beta_1}{\beta}
(\beta , \gamma )$
so
\begin{equation}
A_{21}=\epsilon A_{12},\;\;\epsilon=\frac{\beta_1}{\beta},
\label{s4.3}
\end{equation}
where $\epsilon\neq 0$ as $A_{21}\neq 0$. In the same way, if
$\beta_1\neq 0$, then we get the following  relation
\begin{equation}
A_{12}=\epsilon' A_{21},\;\;\epsilon'=\frac{\beta}{\beta_1}\neq 0.
\end{equation}
Thus, in both cases we have (\ref{s4.3}) with $\epsilon\neq 0$ and $\beta$, 
$\beta_1\neq 0$.\

If $\alpha$, $\gamma\neq 0$ then the conjugation by a diagonal matrix
$diag(1,\alpha,\alpha\beta,\alpha\beta\gamma)$ does not change
$A_{11}$, while $A_{12}$ is reduced to $e_{12}+e_{23}+e_{34}$; i.e. we
obtain $S4a.$ In this case the algebra $\Re $ contains the elements $e_{11}$,
$e_{22}$, $e_{33}$, $e_{44}$ ( as first four powers of $A_{11}$ are linearly
independent) and $e_{12}$=$e_{11}A_{12}$, $e_{23}$=$e_{22}A_{12}$,
$e_{34}$=$e_{33}A_{13}$ as well as elements $e_{13}$=$e_{12}e_{23}$,
$e_{14}$=$e_{13}e_{34}$, $e_{24}$=$e_{23}e_{34}$. This means that $\Re $ 
contains all triangular $4\times 4$ matrices and 
it has the maximal possible dimension.

If $\alpha\neq 0$, $\gamma$=$0$ then, the conjugation by a diagonal matrix 
$diag(1,\alpha,\alpha\beta,1),$ gives us the representation $S4b.$
The algebra
$\Re $ is generated by the elements $e_{11}$, $e_{22}$, $e_{33}$, $e_{44}$,
$e_{12}$, $e_{23}$, $e_{13}$. So this is a 7-dimensional algebra 
of matrices  (see the Table)
which is isomorphic to a direct sum $T_3\oplus C$. The centralizer
consists of diagonal matrices $diag(\beta,\beta,\beta,\delta)$ and is
isomorphic to $C\oplus C.$ 
Formulas (\ref{matrix}) show that the action defined by this 
representation have only one basic invariant; i.e.
\begin{equation}
I=(1-\gamma_0)(1-i\gamma_{12}).
\end{equation}

If $\gamma\neq 0$, $\alpha=0,$ then the conjugation by a diagonal matrix
$diag(1,1,\beta,\beta\gamma)$ gives us the representation $S5,$ 
where the parameter $\alpha $ equals $q^3$ (of course the parameter
$\alpha $ in the Table is not the same as in (\ref{s4.1}) and (\ref{s4.2}), 
which is now zero).

The algebra
$\Re $ is generated by $e_{11}$, $e_{22}$, $e_{33}$, $e_{44}$,
$e_{34}$, $e_{23}$, $e_{24}$. This is also a 7-dimensional 
algebra (see the Table)
which is isomorphic to $T_3\oplus C$, with the centralizer 
$\{diag(\gamma,\delta,\delta,\delta)\}$ isomorphic to $C\oplus C$.
By formulae
(\ref{matrix}) we have that the action defined by this 
representation, has only one basic invariant; i.e.
\begin{equation}
I=(1+\gamma_0)(1+i\gamma_{12}).
\end{equation}

{\bf 5. $\bf A_{11}=diag(\alpha, q^2, q, 1)$, 
$\alpha\neq 0, q^{-1}, 1, q, q^2, q^3$.}
In this case by theorem \ref{t2} we have
\begin{equation}
A_{12}=\beta e_{23}+\gamma e_{34},\ \ \ 
A_{21}=\beta_1e_{23}+\gamma_1e_{34}
\end{equation}
and
\begin{equation}
A_{12}A_{21}=\beta\gamma_1e_{24}=A_{21}A_{12}=\beta_1\gamma e_{24}\neq 0.
\label{s5.1}
\end{equation}
So
\begin{equation}
0=\beta\gamma_1-\beta_1\gamma=det
\left(
\begin{array}{cc}
\beta & \gamma \\
\beta_1 & \gamma_1
\end{array}
\right),\ \ 
\beta,\gamma,\beta_1,\gamma_1\neq 0.
\label{s5.2}
\end{equation}
Therefore
\begin{equation}
A_{21}=\epsilon A_{12},\ \ 
\epsilon=\frac{\beta_1}{\beta}\neq 0.
\label{s5.3}
\end{equation}
Now the conjugation by a diagonal matrix $diag(1,1,\beta,\beta\gamma)$ gives us
the representation $S5$ and the action defined by 
this representation
also has only one basic invariant; i.e.
\begin{equation}
I=(1+\gamma_0)(1+i\gamma_{12}).
\end{equation}
{\bf 6. $\bf A_{11}=diag(q^2, q^2, q,1)+e_{12}$.} By theorem \ref{t2},
$A_{12}$=$\alpha e_{13}+\beta e_{34}$,
$A_{21}$=$\alpha_1e_{13}+\beta_1e_{34},$ thus
\begin{equation} 
A_{12}A_{21}=
\alpha\beta_1e_{14}=\alpha_1\beta e_{14}\neq 0,
\label{s6.1}
\end{equation}
which implies that
\begin{equation}
A_{21}=\epsilon A_{12},\ \  \epsilon=\frac{\alpha_1}{\alpha}\neq 0.
\label{s6.2}
\end{equation}
Conjugation by a diagonal matrix  $diag(1,1,\alpha,\alpha\beta)$ does not 
change the matrix $A_{11}$ while $A_{12}$ is reduced by this to 
$e_{13}+e_{34}$ and we have the representation $S6.$\

For the calculation of $\Re $ let us note that
\begin{equation}
A^k_{11}=diag(q^{2k}, q^{2k}, q^k, 1)+2^{k-1}q^{2k-2}e_{12}.
\label{s6.3}
\end{equation}
Evidently a subalgebra generated by $A_{11}$ is contained in a four 
dimensional algebra of matrices generated by $e_{11}+e_{22}$, $e_{33}$,
$e_{44}$, $e_{12}$. In order to show that these algebras are equal to
each other 
it is enough to show that $E$, $A_{11}$, $A^2_{11}$, $A^3_{11}$ are
linearly independent. If this is not the case then $A_{11}$ is a root of some 
polynomial $f(x)$ of degree three. This polynomial must be a divisor of the
characteristic polynomial of $A_{11}$, which is equal to
\begin{equation}
(x-q^2)^2(x-q)(x-1).
\label{s6.4}
\end{equation}
The polynomial (\ref{s6.4}) has just three divisors of degree three and none of them has
$A_{11}$ as its root (if $q\neq \pm 1$).\

Thus the algebra $\Re $ contains the elements $e_{11}+e_{22}$, $e_{33}$, 
$e_{44}$, $e_{12}$, $e_{13}=A_{12}e_{33}$, $e_{34}$=$A_{12}e_{44}$,
$e_{14}$=$e_{13}e_{34}$ and is generated by these elements 
like a linear space; i.e.
it is the 7-dimensional algebra of matrices presented in the Table. 
Its centralizer is the two 
dimensional algebra $C+Ce_{12}\cong T'_2$ i.e. by formulae (\ref{matrix})
the action defined by this representation has only one basic invariant; i.e.
\begin{equation}
I=(1+\gamma_0)(\gamma_1+i\gamma_2)\gamma_3.
\end{equation}
{\bf 7. $\bf A_{11}=diag(q^2, q, 1, 1)+e_{34}$.} By theorem \ref{t2},
we have $A_{12}$=
$\alpha e_{12}+\beta e_{24}$, $A_{21}=\alpha_1e_{12}+\beta_1e_{24}$
and
\begin{equation}
A_{12}A_{21}=\alpha\beta_1e_{14}=A_{21}A_{12}=\alpha_1\beta e_{14}\neq 0,
\label{s7.1}
\end{equation}
which immediately implies that
\begin{equation}
A_{21}=\epsilon A_{12},\;\;
\epsilon=\frac{\alpha_1}{\alpha}\neq 0.
\label{s7.2}
\end{equation}
A conjugation by the matrix $diag(\alpha^{-1}\beta^{-1}, \beta^{-1}, 1, 1)$
gives the representation $S7.$\

As in the previous case, the algebra $\Re $ is generated by elements
$e_{11}$, $e_{22}$, $e_{33}+e_{44}$, $e_{34}$, $e_{12}$, $e_{24}$, $e_{14},$
i.e. this is the algebra of matrices given in the Table
and its centralizer is equal to $C+e_{34}C\cong T'_2$; i.e. 
the action defined 
by this representation has only one basic invariant
\begin{equation}
I=(1-\gamma_0)(\gamma_1+i\gamma_2)\gamma_3.
\end{equation}
Thus the Theorem \ref{slq} is proved. \hfill$\Box$ 
\section{Actions defined by $GL_q(2,C)$-representations}
   Now let us consider representations of 
$GL_q(2,C)$
which are not equivalent to representations of
$SL_q(2,C).$
If
\begin{equation}
a_{ij}\rightarrow A_{ij}
\label{given}
\end{equation}
is such a representation then we can define a representation of 
$SL_q(2,C)$
setting
$a_{ij}\rightarrow A_{ij}$
for
$i\neq 2$
or
$j\neq 2$
and (see formula (\ref{s}))
\begin{equation}
a_{22}\rightarrow A_{22}^{\prime }=A_{11}^{-1}(1+qA_{12}A_{21}). 
\label{119}
\end{equation}
We will call these two representation {\it connected} to each other. 
Respectively, inner actions defined by connected representations
will also be called {\it connected}.

If we denote by
$D$
the quantum determinant of the given $GL_q(2,C)$-representation
\begin{equation}D=A_{11}A_{22}-qA_{12}A_{21},  \label{120}\end{equation}
then we obtain
\begin{equation}A_{22}=A_{11}^{-1}(D+qA_{12}A_{21})=
A_{22}^{\prime }+A_{11}^{-1}(D-1). \label{121}\end{equation}
In the last formula, 
$D$
is an invertible matrix which commute with all 
$A_{ij}.$
In particular it is an invariant for the connected 
$SL_q(2,C)$-action (\ref{119}).

Inversely, suppose that
\begin{equation}a_{ij}\rightarrow A_{ij}, \ \ i\neq 2 \ {\rm or}\ j\neq 2,\ \ 
a_{22}\rightarrow A_{22}^{\prime }   \label{122}\end{equation}
is a representation of
$SL_q(2,C)$
and 
$D$
is an invertible invariant for the action defined by this representation.
In this case, formula (\ref{121}) defines a representation of
$GL_q(2,C)$
connected with (\ref{122}). 

Thus, every $GL_q(2,C)$-representation is completely defined by 
the connected $SL_q(2,C)$-representation (which can be taken from
the Table) and by an invariant $D$ of the inner action corresponding to
this connected $SL_q(2,C)$-representation (which also can be found in 
the Table).

If $u$
is an invertible invariant for the connected action (\ref{119}), then 
$u$
commutes with all
$A_{ij}, \ (i\neq 2$
or
$j\neq 2).$
Therefore the representation (\ref{given}), where $A_{22}$ is replaced by
\begin{equation}
A_{22}=A^{\prime }_{22}+A_{11}^{-1}(uDu^{-1}-1)
\label{given1}
\end{equation}
is equivalent to (\ref{given}) (see formula (\ref{121}) and
Theorem \ref{t4}. 
In the same way if we multiply
$A_{12}$
and
$D$
by a nonzero scalar 
$\alpha _{2}$
then we will get an equivalent representation with (probably) another 
connected 
$SL_q(2,C)$-representation
and a new
$D$
proportional to the old one.

These considerations show that for the classification, up to the equivalence 
of all representations of
$GL_q(2,C)$ connected with a given $SL_q(2,C)$-representation,
 it is enough to take just one element in every
projective class  of conjugated elements of the group of invertible 
invariants for a connected inner $SL_q(2,C)$-action. 
In particular, if an inner
$SL_q(2,C)$-action
has no nontrivial invariants (like $S2a, S2b, S4a$)
then all 
$GL_q(2,C)$-actions
connected with this action are  
$SL_q(2,C)$-actions.

By Theorem \ref{slq} (see representations marked by $S$ in the Table)
for the algebra of invariants of an $SL_q(2,C)$-action we 
have just three nontrivial possibilities:
$$Inv\cong \left(\matrix{*&*\cr
                 0&*\cr}\right)=T_2,\ \ \
Inv\cong C\oplus C, \ \ \ Inv\cong
           \left(\matrix{\epsilon &*\cr
                         0&\epsilon \cr}\right)=T_{2}^{\prime }.$$

It is easy to see that if a triangular $2\times 2$ matrix has different
nonzero elements on the diagonal, then this matrix is conjugated in $T_2$
with a diagonal one. If this matrix has on its diagonal elements equal to 
$\epsilon $ then this matrix is conjugated
in $T_2$ with a triangular matrix whose nonzero entries 
are all equal to $\epsilon .$ Therefore,
in the first of these three cases we have just two types of 
 possible values for $D,$ which belong to different projective 
classes of conjugated elements:
\begin{equation}D=diag (1, \beta),\ \beta \neq 0,1;\ \ D=\left(\matrix{1&1\cr
                                             0&1\cr}\right). \label{123}\end{equation}

In the second and third cases the algebra $Inv$ is commutative;
 therefore  $D$ can be choosen, respectively, in the forms
\begin{equation}D=1 \oplus \beta ,\ \ \beta \neq 0,1  \label{124}\end{equation}
and
\begin{equation}D=\left(\matrix{1&\xi \cr 0&1\cr}\right),
\ \ \xi \neq 0. \label{125}\end{equation}

Thus, using Theorem \ref{slq} and formula (\ref{121}) with $D$ 
defined by (\ref{123}), (\ref{124}), (\ref{125})   and making 
usual calculations for finding the algebras $\Re $ and their centralizers
one can obtain the following result.

\begin{theorem}. Every representation 
$a_{ij}\rightarrow A_{ij}$
of the algebra
$GL_q(2,C), q^m\neq 1,$ in $C(1,3)$ with nontrivial perturbation
is equivalent to one of the representations given in the Table, 
which contains the operator algebras, the quantum determinants 
and the invariants of
the corresponding inner actions.
\label{glq}
\end{theorem}

In order to show that the Table presents a complete classification it is
necessary to prove that different representations in the Table define
different inner actions.

\begin{theorem}. Two representations given in the Table are equivalent if
and only if they are equal to each other.
\label{=}
\end{theorem}

{\it Proof.} Let
$a_{ij}\rightarrow A_{ij}$ and $a_{ij}\rightarrow A_{ij}^{\prime }$
be equivalent representations given in the Table. By Theorem \ref{t2} 
we have
$$A_{11}^{\prime }=uA_{11}u^{-1}\alpha _1,\ 
  A_{12}^{\prime }=uA_{12}u^{-1}\alpha _2,$$\  
\begin{equation}A_{21}^{\prime }=uA_{21}u^{-1}\alpha _1,\ 
  A_{22}^{\prime }=uA_{22}u^{-1}\alpha _2.   \label{126}\end{equation}
The first of these equations shows that matrices
$A_{11}^{\prime }\alpha_1^{-1}$
and
$A_{11}$
have the same Jordan Normal Form. In particular, if one of the matrices
$A_{11}, A_{11}^{\prime }$
is diagonal then so is the other one and
$A_{11}^{\prime }\alpha _1^{-1}$
can be obtained from 
$A_{11}$
by permutation of its diagonal elements
(note that all matrices $A_{11}$ in the Table have a Jordan Normal Form). 
It easy to see that no one pair of 
{\it different} diagonal matrices $A_{11}$ from the Table satisfies 
this property (here it is essential that $q^4, q^3\neq 1,$
and
$\alpha \neq q^{-1}$
in $S5$). Thus, in these case
$A_{11}^{\prime }=A_{11}.$

If both matrices
$A_{11}, A_{11}^{\prime }$
are not diagonal and
$A_{11}\neq A_{11}^{\prime }$
then one of this matrices appears in S6 while another in S7. Let us assume
$A_{11}^{\prime }=diag(q^2,q,1,1)+e_{34}$
and
$A_{11}=diag(q^2,q^2,q,1)+e_{12}$.
The matrix 
$A_{11}$
has a Jordan Normal Form while the Jordan Normal Form of
$A_{11}^{\prime }\alpha _1^{-1}$
is equal to
\begin{equation}diag(\alpha _1^{-1},\alpha _1^{-1}, q^2\alpha _1^{-1},q\alpha _1^{-1})
+e_{12}.  \label{127}\end{equation}
Therefore
$\alpha _1^{-1}=q^2$
and either
$q^2\alpha _1^{-1}=q, q\alpha _1^{-1}=1$
or
$q^2\alpha _1^{-1}=1, q\alpha _1^{-1}=q.$
Both cases are impossible since
$q^4, q^3\neq 1.$
Thus, we have proved that
$A_{11}^{\prime }=A_{11}$
in all cases.

Now we have
$uA_{11}u^{-1}\alpha _1=A_{11}$ or
$uA_{11}\alpha _1=A_{11}u,$ i.e.
$x\rightarrow A_{11}, y\rightarrow u$
is a representation of the $\alpha _1$-spinor
with invertible both terms. This implies that $\alpha _1$
is a root of unity,
$\alpha _1^m=1.$
If
$(\beta _1, \beta _2, \beta _3, \beta _4)$
is a quadruple of eigenvalues (or, more exactly, the main diagonal
of the Jordan Normal Form) of
$A_{11},$
then the eigenvalues of
$A_{11}\alpha _1$
form a quadruple
$(\beta _1\alpha _1, \beta _2\alpha _1, \beta _3\alpha _1, \beta _4\alpha _1)$
and the equality
$u^{-1}A_{11}u=A_{11}\alpha _1$
shows that these two quadruples coincide up to a permutation
\begin{equation}(\beta _1, \beta _2, \beta _3, \beta _4)=
 (\beta _{\pi (1)}\alpha _1, \beta _{\pi (2)}\alpha _1, 
  \beta _{\pi (3)}\alpha _1, \beta _{\pi (4)}\alpha _1),  \label{128}\end{equation}
where $\pi $ is a permutation of four indices.

Now it is easy to see that no one of the matrices
$A_{11}$
given in the Table satisfies (\ref{128}) for 
$\alpha_1^m=1, \alpha_1\neq 1.$
Thus we have shown that
$\alpha _1=1$
and $u$ commutes with $A_{11}.$

Let us prove that 
$\alpha _2=1.$
The following relation is valid for quantum determinants
\begin{equation}D^{\prime }=uDu^{-1}\alpha _1\alpha _2=uDu^{-1}\alpha _2.  \label{129}\end{equation}
This relation shows that quadruples of eigenvalues of
$D^{\prime }$ and $D$ are connected by the following relation
\begin{equation}(\beta _1^{\prime }, \beta _2^{\prime }, \beta _3^{\prime }, 
\beta _4^{\prime })=
 (\beta _{\pi (1)}\alpha _2, \beta _{\pi (2)}\alpha _2, 
  \beta _{\pi (3)}\alpha _2, \beta _{\pi (4)}\alpha _2).  \label{130}\end{equation}
It is easy to see that for every quantum determinant from the Table
either all four eigenvalues are equal to 1 or three of them are equal to 1
while the fourth one is arbitrary (like in $G1a, G2, G3a, G4, G5$).
In both cases the relations (\ref{130}) are possible only if 
$\alpha _2=1.$

Assume now that no one of our two representations belongs to
$S2a, S2a^\prime , S2b.$
 
All the representations from the Table not belonging to these groups
have
$A_{12}$
with no parameters. If
$A_{11}\neq diag(q^3, q^2, q, 1)$
then this fact immediately implies that $A_{12}^{\prime }=A_{12}.$
For 
$A_{11}=diag(q^3, q^2, q, 1)$ 
the element
$u$
commutes with this diagonal matrix, therefore it is itself 
a diagonal matrix. We see that no different values for
$A_{12}$
(that is 
$e_{12}+e_{23}+e_{34}; \ e_{12}+e_{23}; \ e_{23}+e_{34}$)
are conjugated by a diagonal matrix, so again 
$A_{12}^{\prime }=A_{12}$
and in all the cases 
$uA_{12}u^{-1}=A_{12},$
i.e. $u$ commutes both with 
$A_{11}$
and 
$A_{12}.$
In all of the cases under consideraton
$A_{21}=\alpha A_{12},$
therefore
\begin{equation}
A_{21}^{\prime }=uA_{21}u^{-1}=\alpha uA_{12}u^{-1}=\alpha A_{12}=A_{21}.
\end{equation}
By formula (\ref{119}) this implies that connected
$SL_q(2,C)$-representations
coincide and
$u$
belongs to the algebra of invariants of the connected
$SL_q(2,C)$-representation. So the relation
$D^{\prime }=uDu^{-1}$
implies
$D^{\prime }=D$ 
as by the choice of
$D$ 
(see (\ref{123}), (\ref{124}), (\ref{125}))
different 
$D, D^{\prime }$ are not conjugated in the algebra of invariants of
the connected 
$SL_q(2,C)$-representation. 
Now, formula (\ref{121}) shows that
$A_{22}^{\prime }=A_{22}$
and the representations coincide.

Finally, we have to consider representations with
$A_{11}=diag(q^2, q, q, 1).$
Representations from distinct groups $S2a, S2a^{\prime },\; S2b,\; 
S2b^{\prime },\; G2$ have nonisomorphic algebras 
$\Re $
and therefore they cannot be equivalent (recall that by (\ref{126}) 
the algebras
$\Re $
for equivalent representations are conjugated by $u$). Thus, our 
representations
belong to the same group and we need to consider the first three groups
(in $S2b^{\prime }$ and $G2$ matrices 
$A_{12}$
have no parameters and
$A_{11}, A_{12}$
generate
$\Re $
for the connected
$SL_q(2,C)$-representation).

{$\bf S2a, S2a^\prime $.} 
In these cases $A_{22}$ has no parameters, therefore
$A_{22}^{\prime }=A_{22}$
and the matrix $u$ commutes with 
$A_{11}, A_{22}$.
This means that
$u=diag(1, M, 1),$
where $M$ is an invertible
$2\times 2$ matrix.

If $A_{12}=\alpha e_{12}+e_{13}+e_{24}$
and
$A_{12}^{\prime }=\alpha ^{\prime }e_{12}+e_{13}+e_{24},$
then the relation
$uA_{12}=A_{12}^{\prime }u$
implies
\begin{equation}M\left(\matrix{1\cr 0\cr}\right)=\left(\matrix{1\cr 0\cr}
\right), \ \ \ 
   (\alpha ,1)=(\alpha ^{\prime },1)M.  \label{131}\end{equation}
The first of these equations shows that
$M=\left(\matrix{1&*\cr 0&*\cr }\right)$ 
and in this case
$(\alpha ^{\prime }, 1)M=(\alpha ^{\prime }, *),$
so
$\alpha ^{\prime }=\alpha $
by the second equation of (\ref{131}).

In the same way, if
$A_{21}=e_{12}+\beta e_{13}+e_{34}$
and
$A_{21}^{\prime }=e_{12}+\beta ^{\prime }e_{13}+e_{34},$
then
\begin{equation}(1, \beta )M=(1, \beta ^{\prime }), \ \ \ 
  M\left(\matrix{0\cr 1\cr}\right)=
\left(\matrix{0\cr 1\cr}\right)  \label{132}\end{equation}
and
\begin{equation}\beta ^{\prime }=(1, \beta ^{\prime })\left(\matrix{0\cr 1\cr}\right)=
(1, \beta )M\left(\matrix{0\cr 1\cr}\right)=
(1, \beta )\left(\matrix{0\cr 1\cr}\right)=\beta . \label{133}\end{equation}

{$\bf S2b.$} In this case
$A_{12}$
has no parameters, so
$A_{12}^{\prime }=A_{12}$
and $u$ commutes with 
$A_{11}, A_{12}$.
An algebra generated by
$A_{11}, A_{12}$
of $S2b$ is equal to
$\Re $ 
for $S2b^{\prime }$, thus $u$ is an invariant for $S2b^{\prime }$ i.e.
$u=diag(\beta , \beta , \delta , \beta )$ (see the Table).
Conjugation by $u$ of
$A_{21}=\alpha e_{12}+\alpha e_{24}+e_{13}$
gives
$A_{21}^{\prime }=uA_{12}u^{-1}=\alpha e_{12}+\alpha e_{24}+
         \alpha \delta ^{-1}e_{13}$,
while
$A_{21}^{\prime }=\alpha ^{\prime }e_{12}+\alpha ^{\prime }e_{24}+e_{13}$,
therefore
$\alpha \delta ^{-1}=1$
and
$\alpha =\alpha ^{\prime }$.        \hfill$\Box$

{\it NOTE}.\  In our proofs we do not need the 
fact that every finite dimensional rerpresentation of
$GL_q(2,C)$ is triangular. We need this fact only for four 
or less dimensional
representations. For this we need only restrictions
$q^6, q^8\neq 1$ and $q$ can be a root of unity. 
Our classification is correct for $q^6, q^8\neq 1$
(for example if $q=exp({2\pi i\over 5})$).
It is easy to see that this classification is not valid if
$q=\pm i$ or $q=\pm exp({2\pi i\over 3})$
as in this cases there exist, respectively, four and three dimensional 
irreducible representations of
$SL_q(2,C)$ with $\dim \Re =16, 9,$ respectively. Nevertheless we do 
not know if this classification is correct for
$q=\pm exp({\pi i\over 3})$ or $q^4=-1$.
\section{Table and corollaries} 

In the Table we  have denoted by $S$ (followed by some symbols) the
 representations of $SL_q(2,C)$ and by $G$ (followed by the same symbols) 
representations of $GL_q(2,C),$ 
connected with the corresponding representation of $SL_q(2,C).$ 
We have shown in the Table five ingredients for every representation:
the values of $A_{ij},$ the matrix form of the operator algebra 
$\Re ,$ its dimension,
the invariants of the inner action defined by this representation and 
the value of
the quantum determinant. If the value of some $A_{ij}$ is not shown 
for a particular case then this 
means that it  coincides with the value of the previous 
representation in the Table. By bold style
we have denote seven representations which start groups with the same values
of $A_{11}.$ 
\newpage

\ 

\

\

\

\centerline{\bf Table }

\            

\vbox{\offinterlineskip
\hrule
\halign{&\vrule#&\strut \quad \hfil# \hfil
 &\vrule#&\quad \hfil# \hfil&\vrule#& \quad# \hfil \cr
height5pt&\omit&& \omit&& \omit&\cr&
        {\bf S1}   
&&
$\begin{array}{l}
        \bf A_{11}=diag(q^2,q,1,1)\\ A_{12}=e_{12}+e_{23}\\
A_{21}=\alpha A_{12}, \ \alpha \neq 0\\ 
A_{22}=A_{11}^{-1}+\alpha q^{-1}e_{13} \end{array}$
&& 
        $\Re =\left(\matrix{*&*&*&0\cr 0&*&*&0\cr 0&0&
\epsilon &0\cr 0&0&0&\epsilon \cr}
\right)\cong T_3$&\cr
height5pt&\omit&& \omit&&\omit&\cr
\noalign{\hrule }
height5pt&\omit&& \omit&&\omit&\cr    
&$\matrix{dim \Re \cr 6 \cr }$ &&
         $\begin{array}{l}
{\rm Invariants}\cong T_2\\
         diag(\beta ,\beta ,\beta , \gamma)+\delta e_{43}\\
I_1=(1-\gamma _0)(-\gamma _1+i\gamma _2)\gamma _3\\
I_2=(1-\gamma _0)(1-i\gamma_{12})\end{array}$  && 
         $det_q=1$ &\cr
height5pt&\omit&& \omit&& \omit&\cr
\noalign{\hrule }
height2pt&\omit&& \omit&& \omit&\cr
\noalign{\hrule }
height5pt&\omit&& \omit&& \omit&\cr&
         G1a   
&&
$\begin{array}{c}
A_{22}=A_{11}^{-1}+\alpha q^{-1}e_{13}+\beta e_{44}\\
\alpha ,\beta \neq 0,\ \  \beta \neq -1 \end{array}
$
&& 
        $\Re =\left(\matrix{*&*&*&0\cr 0&*&*&0\cr 0&0&
* &0\cr 0&0&0&* \cr}
\right)\cong T_3\oplus C$&\cr
height5pt&\omit&& \omit&&\omit&\cr
\noalign{\hrule }
height5pt&\omit&& \omit&&\omit&\cr    
&
$\matrix{dim \Re \cr 7 \cr }$
 &&$\begin{array}{l}
     {\rm Invariants}\cong C\oplus C\\
         diag(\gamma ,\gamma ,\gamma ,\delta )\\
I=(1-\gamma _0)(1-i\gamma _{12})\end{array}$  && 
         $det_q=1+\beta e_{44}$ &\cr
height5pt&\omit&& \omit&& \omit&\cr
\noalign{\hrule }
height2pt&\omit&& \omit&& \omit&\cr
\noalign{\hrule }
height5pt&\omit&& \omit&& \omit&\cr&
G1b   
&&
$\begin{array}{c}
         A_{22}=A_{11}^{-1}+\alpha q^{-1}e_{13}+e_{43}\\
\alpha \neq 0 \end{array}$

&& 
        $\Re =\left(\matrix{*&*&*&0\cr 0&*&*&0\cr 0&0&
\epsilon &0\cr 0&0&*&\epsilon \cr}
\right)$&\cr
height5pt&\omit&& \omit&&\omit&\cr
\noalign{\hrule }
height5pt&\omit&& \omit&&\omit&\cr    
&$\matrix{dim \Re \cr 7 \cr }$ 
         &&$\begin{array}{l} 
{\rm Invariants}\cong T_2^{\prime }\\
         diag(\beta ,\beta ,\beta ,\beta )+\delta e_{43}\\
I=(1-\gamma _0)(-\gamma _1+i\gamma _2)\gamma _3\end{array}$  && 
         $det_q=1+e_{43}$ &\cr
height5pt&\omit&& \omit&& \omit&\cr
\noalign{\hrule }
height2pt&\omit&& \omit&& \omit&\cr
\noalign{\hrule }}}
\vbox{\offinterlineskip
\hrule
\halign{&\vrule#&\strut \quad \hfil# \hfil
 &\vrule#&\quad \hfil# \hfil&\vrule#& \quad# \hfil \cr
height5pt&\omit&& \omit&& \omit&\cr
        &{\bf S2}a   
&&
$\begin{array}{l}
        \bf A_{11}=diag(q^2,q,q,1)\\ A_{12}=\alpha e_{12}+e_{13}+e_{24}\\
A_{21}=e_{12}+\beta e_{13}+e_{34}\\ 
A_{22}=A_{11}^{-1}+ q^{-1}e_{14}\\
\hfil \alpha \beta \neq 1 \hfil \end{array}
$
&& 
        $\Re =\left(\matrix{*&*&*&*\cr 0&\epsilon &0&*\cr 0&0&
\epsilon &*\cr 0&0&0&* \cr}
\right)$&\cr
height5pt&\omit&& \omit&&\omit&\cr
\noalign{\hrule }
height5pt&\omit&& \omit&&\omit&\cr    
&$\matrix{dim \Re \cr 8 \cr }$ && Invariants$\cong C$
           && 
         $det_q=1$ &\cr
height5pt&\omit&& \omit&& \omit&\cr
\noalign{\hrule }
height2pt&\omit&& \omit&& \omit&\cr
\noalign{\hrule }
height5pt&\omit&& \omit&& \omit&\cr
        & S2a$^\prime $   
&&
$\alpha \beta =1$
&& 
        $\Re =\left(\matrix{*&\alpha \gamma &\gamma &*\cr 0
&\epsilon &0&*\cr 0&0&
\epsilon &*\cr 0&0&0&* \cr}
\right)$&\cr
height5pt&\omit&& \omit&&\omit&\cr
\noalign{\hrule }
height5pt&\omit&& \omit&&\omit&\cr    
&$\matrix{dim \Re \cr 7 \cr }$ && Invariants$\cong C$
           && 
         $det_q=1$ &\cr
height5pt&\omit&& \omit&& \omit&\cr
\noalign{\hrule }
height2pt&\omit&& \omit&& \omit&\cr
\noalign{\hrule }
height5pt&\omit&& \omit&& \omit&\cr
        & S2b   
&&
$\begin{array}{l}
        A_{11}=diag(q^2,q,q,1)\\ A_{12}=e_{12}+e_{24}\\
A_{21}=\alpha e_{12}+\alpha e_{24}+e_{13} \\ 
A_{22}=A_{11}^{-1}+ q^{-1}\alpha e_{14}\\
\hfil \alpha \neq 0 \hfil \end{array}
$
&& 
        $\Re =\left(\matrix{*&*&*&*\cr 0&\epsilon &0&*\cr 0&0&
\epsilon &0\cr 0&0&0&* \cr}
\right)$&\cr
height5pt&\omit&& \omit&&\omit&\cr
\noalign{\hrule }
height5pt&\omit&& \omit&&\omit&\cr    
&$\matrix{dim \Re \cr 7 \cr }$ && Invariants$\cong C$
           && 
         $det_q=1$ &\cr
height5pt&\omit&& \omit&& \omit&\cr
\noalign{\hrule }
height2pt&\omit&& \omit&& \omit&\cr
\noalign{\hrule }
height5pt&\omit&& \omit&& \omit&\cr
        & S2b$^\prime $   
&&
$A_{21}=\alpha e_{12}+\alpha e_{24},\ \  \alpha \neq 0 
$
&& 
        $\Re =\left(\matrix{*&*&0&*\cr 0&\epsilon &0&*\cr 0&0&
\epsilon &0\cr 0&0&0&* \cr}
\right)\cong T_3$&\cr
height5pt&\omit&& \omit&&\omit&\cr
\noalign{\hrule }
height5pt&\omit&& \omit&&\omit&\cr    
&$\matrix{dim \Re \cr 6 \cr }$ && 
$\begin{array}{l}
{\rm Invariants}\cong C\oplus C\\
                            diag(\beta, \beta ,\delta ,\beta )\\
                              I=(1-\gamma _0)(1+i\gamma _{12})\end{array}$
           && 
         $det_q=1$ &\cr
height5pt&\omit&& \omit&& \omit&\cr
\noalign{\hrule }}}
\vbox{\offinterlineskip
\hrule
\halign{&\vrule#&\strut \quad \hfil# \hfil
 &\vrule#&\quad \hfil# \hfil&\vrule#& \quad# \hfil \cr
height5pt&\omit&& \omit&& \omit&\cr
& G2b$^{\prime }$   
&&
$\begin{array}{l}
        A_{11}=diag(q^2,q,q,1)\\
        A_{12}=e_{12}+e_{24}\\
        A_{21}=\alpha A_{12}, \ \alpha\neq 0\\ 
A_{22}=A_{11}^{-1}+ q^{-1}\alpha e_{14}+\beta e_{33}\\
\hfil \beta \neq -q^{-1},\  0 \hfil \end{array}
$
&& 
        $\Re =\left(\matrix{*&*&0&*\cr 0&*&0&*\cr 0&0&
*&0\cr 0&0&0&* \cr}
\right)\cong T_3\oplus C$&\cr
height5pt&\omit&& \omit&&\omit&\cr
\noalign{\hrule }
height5pt&\omit&& \omit&&\omit&\cr    
&$\matrix{dim \Re \cr 7 \cr }$ && 
$\begin{array}{l}
{\rm Invariants}\cong C\oplus C\\
                             diag(\delta ,\delta ,\epsilon ,\delta )\\
                             I=(1-\gamma _0)(1+i\gamma _{12})\end{array}$
           && 
         $det_q=1+q\beta e_{33}$ &\cr
height5pt&\omit&& \omit&& \omit&\cr
\noalign{\hrule }
height0.5pt&\omit&& \omit&& \omit&\cr
\noalign{\hrule }
height0.5pt&\omit&& \omit&& \omit&\cr
\noalign{\hrule }
height5pt&\omit&& \omit&& \omit&\cr&
        {\bf S3}   
&&
$\begin{array}{l}
        \bf A_{11}=diag(q^2,q^2,q,1)\\ 
            A_{12}=e_{13}+e_{34}\\
            A_{21}=\alpha A_{12}, \ \alpha \neq 0\\
            A_{22}=A_{11}^{-1}+\alpha q^{-1}e_{14}
\end{array}$
&& 
        $\Re =\left(\matrix{\epsilon &0&*&*\cr 0&\epsilon &0&0\cr 0&0&
 *&*\cr 0&0&0&* \cr}
\right)\cong T_3$&\cr
height5pt&\omit&& \omit&&\omit&\cr
\noalign{\hrule }
height5pt&\omit&& \omit&&\omit&\cr    
&$\matrix{dim \Re \cr 6 \cr }$ &&
$\begin{array}{l}
         {\rm Invariants}\cong T_2\\
         diag(\beta ,\gamma, \beta, \beta )+\delta e_{12}\\
         I_1=(1+\gamma _0)(1-i\gamma _{12})\\
         I_2=(1+\gamma _0)(\gamma_1+i\gamma_2)\gamma_3\end{array}$  && 
         $det_q=1$ &\cr
height5pt&\omit&& \omit&& \omit&\cr
\noalign{\hrule }
height2pt&\omit&& \omit&& \omit&\cr
\noalign{\hrule }
height5pt&\omit&& \omit&& \omit&\cr&
         G3a   
&&
$\begin{array}{c}
       A_{22}=A_{11}^{-1}+\alpha q^{-1}e_{14}+\beta e_{22}\\
       \alpha ,\beta \neq 0,\ \beta \neq -q^{-2}\end{array}$
&& 
        $\Re =\left(\matrix{*&0&*&*\cr 0&*&0&0\cr 0&0&
*&*\cr 0&0&0&* \cr}
\right)\cong T_3\oplus C$&\cr
height5pt&\omit&& \omit&&\omit&\cr
\noalign{\hrule }
height5pt&\omit&& \omit&&\omit&\cr    
&
$\matrix{dim \Re \cr 7 \cr }$
 &&$\begin{array}{l}
          {\rm Invariants}\cong C\oplus C\\
          diag(\gamma ,\delta ,\gamma ,\gamma )\\
          I=(1+\gamma _0)(1-i\gamma _{12})\end{array}$  && 
         $det_q=1+\beta q^2 e_{22}$ &\cr
height5pt&\omit&& \omit&& \omit&\cr
\noalign{\hrule }
height2pt&\omit&& \omit&& \omit&\cr
\noalign{\hrule }
height5pt&\omit&& \omit&& \omit&\cr&
G3b   
&&
$\begin{array}{c}
         A_{22}=A_{11}^{-1}+\alpha q^{-1}e_{14}+e_{12}\\
         \alpha \neq 0\end{array}$

&& 
        $\Re =\left(\matrix{\epsilon &*&*&*\cr 0&\epsilon &0&0\cr 0&0&
*&*\cr 0&0&0&* \cr}
\right)$&\cr
height5pt&\omit&& \omit&&\omit&\cr
\noalign{\hrule }
height5pt&\omit&& \omit&&\omit&\cr    
&$\matrix{dim \Re \cr 7 \cr }$ 
         &&$\begin{array}{l}
            {\rm Invariants}\cong T_2^{\prime }\\
            diag(\beta ,\beta ,\beta ,\beta )+\delta e_{12}\\
            I=(1+\gamma _0)(\gamma _1+i\gamma _2)\gamma _3 \end{array}$  && 
 $det_q=1+e_{12}$ &\cr
height5pt&\omit&& \omit&& \omit&\cr
\noalign{\hrule }}}
\vbox{\offinterlineskip
\hrule
\halign{&\vrule#&\strut \quad \hfil# \hfil
 &\vrule#&\quad \hfil# \hfil&\vrule#& \quad# \hfil \cr
height5pt&\omit&& \omit&& \omit&\cr&
        {\bf S4}a   
&&
$\begin{array}{l} 
      \bf A_{11}=diag(q^3,q^2,q,1)\\
          A_{12}=e_{12}+e_{23}+e_{34} \\
          A_{21}=\alpha A_{12}, \ \alpha \neq 0\hfil \\
          A_{22}=A_{11}^{-1}+\alpha q^{-2}e_{13}+\alpha q^{-1}e_{24} 
\end{array}$
&& 
        $\Re =\left(\matrix{*&*&*&*\cr 0&*&*&*\cr 0&0&
*&*\cr 0&0&0&* \cr}
\right)\cong T_4$&\cr
height5pt&\omit&& \omit&&\omit&\cr
\noalign{\hrule }
height5pt&\omit&& \omit&&\omit&\cr    
&$\matrix{dim \Re \cr 10 \cr }$ && Invariants$\cong C$
           && 
         $det_q=1$ &\cr
height5pt&\omit&& \omit&& \omit&\cr
\noalign{\hrule }
height2pt&\omit&& \omit&& \omit&\cr
\noalign{\hrule }
height5pt&\omit&& \omit&& \omit&\cr&
         S4b   
&&
$\begin{array}{l}
         A_{11}=diag(q^3,q^2,q,1)\cr 
         A_{12}=e_{12}+e_{23}\\
         A_{21}=\alpha A_{12}, \ \alpha \neq 0\\
         A_{22}=A_{11}^{-1}+\alpha q^{-2}e_{13}\end{array}
$
&& 
        $\Re =\left(\matrix{*&*&*&0\cr 0&*&*&0\cr 0&0&
*&0\cr 0&0&0&* \cr}
\right)\cong T_3\oplus C$&\cr
height5pt&\omit&& \omit&&\omit&\cr
\noalign{\hrule }
height5pt&\omit&& \omit&&\omit&\cr    
&$\matrix{dim \Re \cr 7 \cr }$ 
           && $\begin{array}{l}
              {\rm Invariants}\cong C\oplus C\\
              diag(\beta ,\beta ,\beta ,\delta )\\
              I=(1-\gamma _0)(1-i\gamma _{12})\hfil \end{array}$&&
         $det_q=1$ &\cr
height5pt&\omit&& \omit&& \omit&\cr
\noalign{\hrule }
height2pt&\omit&& \omit&& \omit&\cr
\noalign{\hrule }
height5pt&\omit&& \omit&& \omit&\cr&
         G4b   
&&
$\begin{array}{l}
         A_{11}=diag(q^3,q^2,q,1)\cr 
         A_{12}=e_{12}+e_{23}\\
         A_{21}=\alpha A_{12}, \ \alpha \neq 0\\
         A_{22}=A_{11}^{-1}+\alpha q^{-2}e_{13}+\beta e_{44}\\
\hfil \alpha ,\beta \neq 0,\ \beta \neq -1 \hfil \end{array}$
&& 
        $\Re =\left(\matrix{*&*&*&0\cr 0&*&*&0\cr 0&0&
*&0\cr 0&0&0&* \cr}
\right)\cong T_3\oplus C$&\cr
height5pt&\omit&& \omit&&\omit&\cr
\noalign{\hrule }
height5pt&\omit&& \omit&&\omit&\cr    
&
$\matrix{dim \Re \cr 7 \cr }$
 &&$\begin{array}{l}
         {\rm Invariants}\cong C\oplus C\\
         diag(\gamma ,\gamma ,\gamma ,\delta )\\
         I=(1-\gamma _0)(1-i\gamma _{12})\end{array}$  && 
         $det_q=1+\beta e_{44}$ &\cr
height5pt&\omit&& \omit&& \omit&\cr
\noalign{\hrule }}}
\vbox{\offinterlineskip
\hrule
\halign{&\vrule#&\strut \quad \hfil# \hfil
 &\vrule#&\quad \hfil# \hfil&\vrule#& \quad# \hfil \cr
height5pt&\omit&& \omit&& \omit&\cr&
        {\bf S5}   
&&
$\begin{array}{l}
        \bf A_{11}=diag(\alpha ,q^2,q,1)\\
        \hfil \alpha \neq 0,q^{-1},1,q,q^{2}\hfil \\
        A_{12}=e_{23}+e_{34}\\
        A_{21}=\beta A_{12},\ \beta \neq 0\\ 
        A_{22}=A_{11}^{-1}+\beta q^{-1}e_{24}\end{array}$
&& 
        $\Re =\left(\matrix{*&0&0&0\cr 0&*&*&*\cr 0&0&
*&*\cr 0&0&0&* \cr}
\right)\cong T_3\oplus C$&\cr
height5pt&\omit&& \omit&&\omit&\cr
\noalign{\hrule }
height5pt&\omit&& \omit&&\omit&\cr    
&$\matrix{dim \Re \cr 7 \cr }$ && $\begin{array}{l}
          {\rm Invariants}\cong C\oplus C\\
          diag(\gamma ,\delta ,\delta ,\delta )\\
          I=(1+\gamma _0)(1+i\gamma _{12})\end{array}$
           && 
         $det_q=1$ &\cr
height5pt&\omit&& \omit&& \omit&\cr
\noalign{\hrule }
height2pt&\omit&& \omit&& \omit&\cr
\noalign{\hrule }
height5pt&\omit&& \omit&& \omit&\cr&
         G5  
&&
$\begin{array}{c}
        A_{22}=A_{11}^{-1}+q^{-1}\beta e_{24}+\gamma e_{11}\\
           \beta \neq 0, \ \gamma \neq -\alpha ^{-1}, 0\end{array}$
&& 
        $\Re =\left(\matrix{*&0&0&0\cr 0&*&*&*\cr 0&0&
*&*\cr 0&0&0&* \cr}
\right)\cong T_3\oplus C$&\cr
height5pt&\omit&& \omit&&\omit&\cr
\noalign{\hrule }
height5pt&\omit&& \omit&&\omit&\cr    
&$\matrix{dim \Re \cr 7 \cr }$ && $\begin{array}{l}
            {\rm Invariants}\cong C\oplus C\\
            diag(\delta ,\epsilon ,\epsilon ,\epsilon )\\
            I=(1+\gamma _0)(1+i\gamma _{12})\end{array}$
           && 
         $det_q=1+\alpha \gamma e_{11}$ &\cr
height5pt&\omit&& \omit&& \omit&\cr
\noalign{\hrule }
height0.5pt&\omit&& \omit&& \omit&\cr
\noalign{\hrule }
height0.5pt&\omit&& \omit&& \omit&\cr
\noalign{\hrule }
height5pt&\omit&& \omit&& \omit&\cr
        &{\bf S6}   
&&
$\begin{array}{l}
        \bf A_{11}=diag(q^2 ,q^2,q,1)+e_{12}\\
        A_{12}=e_{13}+e_{34}\\
        A_{21}=\alpha A_{12},\ \alpha \neq 0\\ 
        A_{22}=A_{11}^{-1}+\alpha q^{-1}e_{14} \end{array}$
&& 
        $\Re =\left(\matrix{\epsilon &*&*&*\cr 0&\epsilon &0&0\cr 0&0&
*&*\cr 0&0&0&* \cr}
\right)$&\cr
height5pt&\omit&& \omit&&\omit&\cr
\noalign{\hrule }
height5pt&\omit&& \omit&&\omit&\cr    
&$\matrix{dim \Re \cr 7 \cr }$ && 
$\begin{array}{l}
         {\rm Invariants}\cong T_2^{\prime }\\
         diag(\beta ,\beta ,\beta ,\beta )+\delta e_{12}\\
         I=(1+\gamma _0)(\gamma _1+i\gamma _{2})\gamma _3\end{array}$
           && 
         $det_q=1$ &\cr
height5pt&\omit&& \omit&& \omit&\cr
\noalign{\hrule }
height2pt&\omit&& \omit&& \omit&\cr
\noalign{\hrule }
height5pt&\omit&& \omit&&\omit&\cr&
         G6   
&&
$A_{22}=\left(\matrix{q^{-2} & \xi -1 & 0 & \alpha q^{-1}\cr
               0      & q^{-2}      & 0 & 0            \cr
               0      &    0    & q^{-1}& 0            \cr
               0      &    0    &    0  & 1            \cr}\right)$
       
&& 
        $\Re =\left(\matrix{\epsilon &*&*&*\cr 0&\epsilon &0&0\cr 0&0&
*&*\cr 0&0&0&* \cr}
\right)$&\cr
height5pt&\omit&& \omit&&\omit&\cr
\noalign{\hrule }
height5pt&\omit&& \omit&&\omit&\cr    
&$\matrix{dim \Re \cr 7 \cr }$ && 
$\begin{array}{l}
         {\rm Invariants}\cong T_2^{\prime }\\
         diag(\beta ,\beta ,\beta ,\beta )+\delta e_{12}\\
         I=(1+\gamma _0)(\gamma _1+i\gamma _{2})\gamma _3\end{array}$
           && 
         $det_q=1+q^2\xi e_{12}, \ \xi \neq 0$  &\cr
height5pt&\omit&& \omit&& \omit&\cr
\noalign{\hrule }}}
\vbox{\offinterlineskip
\hrule
\halign{&\vrule#&\strut \quad \hfil# \hfil
 &\vrule#&\quad \hfil# \hfil&\vrule#& \quad
  \hfil# \hfil \cr
height5pt&\omit&& \omit&& \omit&\cr&
        {\bf S7}   
&&
$\begin{array}{l}
        \bf A_{11}=diag(q^2,q,1,1)+e_{34}\\ A_{12}=e_{12}+e_{24}\\
A_{21}=\alpha A_{12},\ \alpha\neq 0\\ 
A_{22}=A_{11}^{-1}+\alpha q^{-1}e_{14}\end{array}$
&& 
        $\Re =\left(\matrix{*&*&0&*\cr 0&*&0&*\cr 0&0&
\epsilon &*\cr 0&0&0&\epsilon \cr}
\right)$&\cr
height5pt&\omit&& \omit&&\omit&\cr
\noalign{\hrule }
height5pt&\omit&& \omit&&\omit&\cr    
&$\matrix{dim \Re \cr 7 \cr }$ && 
$\begin{array}{l}
        {\rm Invariants}\cong T_2^{\prime }\\
        diag(\beta ,\beta ,\beta ,\beta )+\delta e_{34}\\
        I=(1-\gamma _0)(\gamma _1+i\gamma _{2})\gamma _3
 \end{array}$
           && 
         $det_q=1$ &\cr
height5pt&\omit&& \omit&& \omit&\cr
\noalign{\hrule }
height2pt&\omit&& \omit&& \omit&\cr
\noalign{\hrule }
height5pt&\omit&& \omit&& \omit&\cr&
         G7   
&&
$A_{22}=\left(\matrix{q^{-2}&0&0&\alpha q^{-1}\cr
                              0&q^{-1}&0&0\cr
                              0&0&1&\xi -1\cr
                              0&0&0&1\cr }\right)$
&& 
        $\Re =\left(\matrix{*&*&0&*\cr 0&*&0&*\cr 0&0&
\epsilon &*\cr 0&0&0&\epsilon \cr}
\right)$&\cr
height5pt&\omit&& \omit&&\omit&\cr
\noalign{\hrule }
height5pt&\omit&& \omit&&\omit&\cr    
&$\matrix{dim \Re \cr 7 \cr }$ && 
$\begin{array}{l} 
       {\rm Invariants}\cong T_2^{\prime }\\
       diag(\beta ,\beta ,\beta ,\beta )+\delta e_{34}\\
       I=(1-\gamma _0)(\gamma _1+i\gamma _{2})\gamma _3\end{array}$
           && 
         $det_q=1+\xi e_{34}, \ \xi \neq 0 $ &\cr
height5pt&\omit&& \omit&& \omit&\cr
\noalign{\hrule }}}

\

\

From the Table, it is easy to see that the following 
result is valid (roughly speaking this means that the quantum determinants
are the only quantum invariants).
\begin{corollary}.
If an inner action defined by a representation of $GL_q(1,C),$ \ $q^m\neq 1,$
 in $C(1,3)$ with nontrivial perturbation is not an action of $SL_q(2,C),$
then it has only one basic invariant which can be taken to be equal
to the quantum determinant.
\label{gli}
\end{corollary} 

Of course for $SL_q(2,C)$-actions this result is not valid 
(in this case the quantum determinant equals 1). Nevertheless,
we have seen (formulae (\ref{122}, \ref{121})) that for every 
invariant $D$ there exists a connected
$GL_q(2,C)$-action defined by a representation for which the quantum 
determinant equals $D$ (even if the initial $SL_q(2,C)$-action is 
defined by a representation with trivial perturbation). 
So we can formulate
\begin{corollary}.
Every invariant of the inner $SL_q(2,C)$-action on $C(1,3)$ is equal to a
quantum determinant of a connected $GL_q(2,C)$-action.
\label{sli}
\end{corollary}
\section{Acknowledgments}
We wish to thank Dr. Zbigniew Oziewicz for helpful discussions.
VKK wishes to thank SNI and CONACYT-M\'exico for its support 
under grant No. 940411-R96 and also the Russian
Foundation for Fundamental Research, grant 95-01-01356. 
SRR wishes to thank
CONACYT-M\'exico for partial support under grant No 4336-E.

\end{document}